\documentclass[11pt]{amsart}

\theoremstyle{plain}
\begin{document}
\newcommand{\RR}{{\mathbb R}}
 \newcommand{\Hi}{{\mathbb H}}
\newcommand{\Ss}{{\mathbb S}}
\newcommand{\N}{{\mathbb N}}
\newcommand{\Rn}{{\RR^n}}
\newcommand{\ieq}{\begin{equation}}
\newcommand{\eeq}{\end{equation}}
\newcommand{\ieqa}{\begin{eqnarray}}
\newcommand{\eeqa}{\end{eqnarray}}
\newcommand{\ieqas}{\begin{eqnarray*}}
\newcommand{\eeqas}{\end{eqnarray*}}
\newcommand{\Bo}{\put(260,0){\rule{2mm}{2mm}}\\}
\def\L#1{\label{#1}} \def\R#1{{\rm (\ref{#1})}}
\newtheorem{defin}{Definition}[section]
\newtheorem{theorem}{Theorem}[section]
    \newtheorem{stat}{Statement}[section]
\newtheorem{lemma}{Lemma}[section]
\newtheorem{remark}{Remark}[section]
\newtheorem{definition}{Definition}[section]
\newtheorem{assumption}{Assumption}[section]
    \newtheorem{appendice}{Appendix}[section]
\newtheorem{proposition}{Proposition}[section]
\newtheorem{corollary}{Corollary}[section]
\newtheorem{step}{Step}[section]
\newcommand{\res}{\mathop{\hbox{\vrule height 7pt width .5pt depth
0pt \vrule height .5pt width 6pt depth 0pt}}\nolimits} \def\at#1{{\bf
#1}: } \def\att#1#2{{\bf #1}, {\bf #2}: } \def\attt#1#2#3{{\bf #1},
{\bf #2}, {\bf #3}: } \def\atttt#1#2#3#4{{\bf #1}, {\bf #2}, {\bf
#3},{\bf #4}: } \def\aug#1#2{\frac{\displaystyle #1}{\displaystyle
#2}} \def\figura#1#2{ \begin{figure}[ht] \vspace{#1} \caption{#2}
\end{figure}} \def\B#1{\bibitem{#1}} \def\q{\int_{\Omega^\sharp}}
\def\z{\int_{B_{\bar{\rho}}}\underline{\nu}\nabla (w+K_{c})\cdot
\nabla h} \def\a{\int_{B_{\bar{\rho}}}} \def\b{\cdot\aug{x}{\|x\|}}
\def\n{\underline{\nu}} \def\d{\int_{B_{r}}}
\def\e{\int_{B_{\rho_{j}}}} \def\LL{{\mathcal L}}
\def\itr{\mathrm{Int}\,}
\def\D{{\mathcal D}}
 \def\tg{\tilde{g}}
\def\A{{\mathcal A}}
\def\S{{\mathcal S}}
\def\H{{\mathcal H}}
\def\M{{\mathcal M}}
\def\T{{\mathcal T}}
\def\U{{\mathcal U}}
\def\N{{\mathcal N}}
\def\I{{\mathcal I}}
\def\F{{\mathcal F}}
\def\J{{\mathcal J}}
\def\E{{\mathcal E}}
\def\F{{\mathcal F}}
\def\G{{\mathcal G}}
\def\HH{{\mathcal H}}
\def\W{{\mathcal W}}
\def\H{\D^{2*}_{X}}
\def\d{d^X_M }
\def\LL{{\mathcal L}}
\def\H{{\mathcal H}}
\def\HH{{\mathcal H}}
\def\itr{\mathrm{Int}\,}
\def\vah{\mbox{var}_\Hi}
\def\vahh{\mbox{var}_\Hi^1}
\def\vax{\mbox{var}_X^1}
\def\va{\mbox{var}}
\def\SS{{\mathcal S}}
 \def\Y{{\mathcal Y}}
\def\length{{l_\Hi}}
\newcommand{\average}{{\mathchoice {\kern1ex\vcenter{\hrule
height.4pt width 6pt depth0pt}
\kern-11pt} {\kern1ex\vcenter{\hrule height.4pt width 4.3pt
depth0pt} \kern-7pt} {} {} }}
\newcommand{\ave}{\average\int}
\title[The Hessian of the
distance from a surface in the Heisenberg group]{The Hessian of the
distance from a surface in the Heisenberg group\footnote{MSC 49Q15, 53C17, 53C22}}
\author{Nicola Arcozzi, Fausto Ferrari}
\date{7/12/2005}
\begin{abstract}
Given a smooth surface $S$ in the Heisenberg group, we 
compute the Hessian of the function measuring the Carnot-Charath\'eodory
distance from $S$ in terms of the Mean Curvature of $S$ and of an ``imaginary
curvature'' which was introduced in \cite{AF1} in order to find the geodesics 
which are metrically normal to $S$. Explicit formulae are given when $S$ is a 
plane or the metric sphere.
\end{abstract}
\maketitle
\tableofcontents

\section{Introduction}
In this article we continue to study the properties of the function ``signed distance
from a surface  $S$'', $\delta_S$,  in the Heisenberg group, started in \cite{AF1}. 
In particular we are interested in the horizontal Hessian of $\delta_S$ and in the related
 notions of curvature for $S$.
Let 
$\Hi=\Hi^1$ be
${\mathbb R}^3$ with the Heisenberg
group structure, and let $d$ be the associated Carnot-Charath\'eodory distance.
See Section \ref{definizioni}  for the definitions. If $S\subset\Hi$
is closed, the distance from a point $P$ to $S$ is
$$
d_S(P)=\inf_{Q\in S}d(P,Q)
$$
Here, we consider the case where $S=\partial\Omega$ is the  $C^3$ boundary, in the Euclidean sense,
of an open subset of $\Hi$. It would be natural to consider  $C^2$ surfaces in
the Euclidean sense, since such surfaces
satisfy a internal/external Heisenberg sphere condition at any
non-characteristic point, see \cite{AF1}. However, our elementary approach requires a degree more 
of regularity. This is a minor nuisance, since here we are interested in geometric properties of 
surfaces.

The {\bf signed distance}
from $S$ is
\begin{equation}
\label{eqsigndist}
\delta_{S}(P)=
\begin{cases}
-d_{S}(P) & {\mbox{if}}\ P\in\Omega\cr
d_{S}(P) & {\mbox{if}}\ P\notin\Omega
\end{cases}
\end{equation}
Observe that, given $S$, the signed distance $\delta_S$ is defined modulo a sign,
which corresponds to the choice of one of the two open sets having
boundary $S$, or, equivalently, to the choice of an orientation for $S$.

It was shown in \cite{MoSe} that, if $S$ is a closed subset of $\Hi$, then $d_S$ is $a.e.$ solution of the 
eikonal equation 
\begin{equation}\label{eikonal}
|\nabla_\Hi d_S|=1,
\end{equation}
where $\nabla_\Hi=(X,Y) $ 
is the horizontal gradient. More on 
this line of investigation is in \cite{BaCa}, \cite{CaSi}, \cite{Dra}, \cite{Si}. 
Recall that the left invariant vector 
fields
$$
X=\partial_x+2y\partial_t,\ Y=\partial_y-2x\partial_t
$$
are a orthonormal basis for the horizontal bundle ${\mathcal H}$ of $\Hi$.
Concerning the properties of the distance function from a point and the distance function 
from a surface,
we recall the paper by Agrachev and Gautier, \cite{AG} and the
contribution by Vershik and Gershkovich, which is surveyed in \cite{VG}. In \cite{AF1}, 
we proved several regularity results for 
the function $\delta_S$, in particular that, if $S$ is $C^2$ in the Euclidean sense, then 
$\nabla_\Hi\delta_S$ is $C^1$ in the Euclidean sense in an open neighborhood of $S\setminus Char(S)$,
where $Char(S)$ is the characteristic set of $S$.

In this article, we move a step forward in the understanding of the higher regularity of $\delta_S$.
The {\bf horizontal Hessian} of a smooth function $f$ defined on an open subset of $\Hi$ is the
matrix
$$
Hess_\Hi f=
\left[\begin{array}
{cc}
XXf & YXf \\
XYf & YYf 
\end{array}\right].
$$
See, for instance, \cite{DGN} . We denote by $I$ and $J$ the identity and the simplectic 
matrix, respectively,
$$
I=\left[\begin{array}
{cc}
1 & 0 \\
0 & 1 
\end{array}\right],
J=\left[\begin{array}
{cc}
0 & 1 \\
-1 & 0 
\end{array}\right].
$$
\begin{theorem}
\label{Tmain}
Let $S=\partial\Omega$ be boundary of a $C^3$ open subset  of $\Hi$ and let $P$ be a 
noncharacteristic point of $S$.
Then,
\begin{equation}
\label{eqhessnsi}
\mbox{Hess}_\Hi\ \delta_S(P)=v_PS\otimes v_PS \cdot (h_S(P)I+p_S(P)J),
\end{equation}
where $v_PS$ is the unit horizontal vector tangent to $S$ at $P$, $h_S(P)$ is the mean curvature 
of $S$
at $P$ and $p_S(P)$ is the curvature of the oriented metric normal to $S$ at $P$.
\end{theorem}
We shall also refer to $h_S(P)$ as the {\bf real curvature} and to $p_S(P)$ as the 
{\bf imaginary curvature} of $S$ at $P$, respectively.
Some explanation about the terminology is in order. Complete definitions will be given in Section 
2.

If $P$ is a noncharacteristic point of $S$, the intersection of the space tangent to $S$ at $P$, 
$T_PS$,
with the horizontal fiber ${\mathcal H}_P$ is a one dimensional linear space $V_PS$, 
the {\bf horizontal direction tangent to $S$ at $P$}. The {\bf unit horizontal vector tangent to 
$S$ at $P$}
is the unit vector spanning $V_PS$, which is defined modulo a sign. There is a rich literature on 
the 
calculus on surfaces in $\Hi$. Just to mention a few, we refer to \cite{Pansu},
\cite{FSSC1}, \cite{FSSC2}, \cite{FSSC3}, \cite{GP}, \cite{Pa1}, \cite{Pa2}.
 In \cite{FSSC2}, in particular, surfaces are 
defined
in terms of graphs of functions, with the aid of an implicit function theorem, an approach which 
is 
especially interesting in that it provides an intrinsic definition of $C^1$ surfaces in the 
Heisenberg 
sense.  
We remark that 
$v_PS$ is a first order differential object in both the Heisenberg and the Euclidean sense.

If the open set $\Omega\subset\Hi$ and $S=\partial\Omega$ have, near $P$, 
the smooth parametrizations $g(z,t)<0$ and $g(z,t)=0$, respectively, then
the analytic expression for $p_S(P)$ is
\begin{equation}
\label{reggio}
p_S(P)=-\frac{[X,Y]g(P)}{|\nabla_\Hi g(P)|}.
\end{equation}
We say that the parametrization $g$ is {\it compatible with the orientation} of $S$.
Observe that $p_S(P)$ is a $C^1$ object in the Euclidean sense, since $[X,Y]g=-4\partial_tg$, 
while 
it is intermediate, in a sense, between $C^1$ and $C^2$ in the Heisenberg sense. It was shown in 
\cite{AF1} how $p_S$ is a quantity peculiarly related to the Heisenberg geometry. Let 
${\mathcal N}^+_PS=\gamma$ 
be the {\bf oriented metric normal to $S$ at $P$}, i.e., the geodesic arc $\gamma$ 
leaving $\Omega$ at $P$ such that, for small $\sigma>0$, $\delta_S(\gamma(\sigma))=\sigma$.
The lifetime $\tau$ of the geodesic $\gamma$,  the maximum amount of time over which $\gamma$
is length-minimizing, is $\tau=\frac{\pi}{|p_S(P)|}$. The sign of $p_S(P)$ is positive or
negative according to the fact that ${\mathcal N}^+_PS$ points ``upward'' (the $t$-coordinate of
$P^{-1}\cdot{\mathcal N}_P^+S(\sigma)$ 
increases with $\sigma$) or ``downward''. $p_S(P)=0$ if ${\mathcal N}_P^+S$ 
is a Euclidean straight line. The number
$p_S(P)$ is a sort of curvature for ${\mathcal N}^+_PS$. Theorem \ref{Tmain} says, in a way, 
that $p_S(P)$ is also a sort of 
curvature of the surface $S$. Observe that the quantity $p_S$ does not appear
in the expression of $\mbox{Hess}^{sym}_\Hi=\mbox{Hess}_\Hi+\mbox{Hess}^{tr}_\Hi$, 
the {\bf symmetrized horizontal Hessian} 
$\mbox{Hess}^{sym}_\Hi\ \delta_S$
 of $\delta_S$,
 the superscript $^{tr}$ denoting 
transposition of matrices. 

The {\bf mean curvature} $h_S(P)$ of a surface $S$ at $P$ is an object much studied in the recent
literature 
\cite{CHMY}, 
\cite{GP}, \cite{Pa1}, \cite{Pa2}, \cite{BC}.  If $g$ is a parametrization of $S$ near the 
noncharacteristic point $P$,  which is compatible with the orientation of $S$,
the mean curvature is defined by
\begin{equation}
h_S(P)=
X\left(\frac{Xg(P)}{|\nabla_\Hi g(P)|}\right)+Y\left(\frac{Yg(P)}{|\nabla_\Hi g(P)|}\right).
\end{equation}
Observe that, if $S=\{(z,t):\ g(z,t)=0\}$, and $\Omega=\{(z,t):\ g(z,t)<0\}$, then
 $$
h_S(P)=div_\Hi N_PS,
$$ where $N_PS=\frac{\nabla_\Hi g(P)}{|\nabla_\Hi g(P)|}$ 
is the intrinsic  outer unit normal to $S$ at $P$.  
Observe that, by Theorem \ref{Tmain}, it is tautological that 
$$
Trace(\mbox{Hess}_\Hi\ \delta_S(P))=h_S(P).
$$
In Section \ref{secmean}, with some applications in perspective, we show that there are other
equivalent definitions for $h_S(P)$. On the one hand, $h_S(P)$ is the eigenvalue of a suitable
{\bf horizontal Weingarten map} on $S$ at $P$, in analogy with the Euclidean case 
(Theorem \ref{Tcurvinv}). 
On the other hand, if $\Gamma$ is the unique horizontal curve on $S$ through $P$ and $\pi(\Gamma)$
is its vertical projection onto the plane $t=0$ in $\Hi$, then $h_S(P)$ is the {\it Euclidean} 
curvature 
of $\pi(\Gamma)$ at $\pi(P)$ (Theorem \ref{Tcurveuc}). Some of these facts have been 
independently proved in
a similar form in \cite{BC}, \cite{GP}, \cite{CHMY}, \cite{CMS}.
In \cite{GP} it is also proved that $h_S$ is the trace of a second order matrix which is different 
from the horizontal Hessian of $\delta_S$. For a more general treatment of the parallel transport 
see, e.g., \cite{J} for the Riemannian case and \cite{M} for the subriemannian one.

As an application of Theorem \ref{Tmain} and Theorem \ref{Tcurveuc}, in Theorem \ref{Thesssigma}
we give the explicit expression for the horizontal Hessian of the function measuring the 
Carnot-Charath\'eodory distance from a point.

There are two reasons why we restricted our analysis to the Heisenberg group. The first is out of
necessity: the proofs in this paper are based on the results in \cite{AF1}, which are stated and 
proved in this setting. The second is that the Heisenberg group provides the simplest nontrivial 
example of Carnot-Charath\'eodory group and it is important in itself, for instance in modelling 
the visual cortex \cite{CS}, \cite{CMS}, \cite{SCM}. It is interesting that these new 
applications require a finer study of surfaces in the Heisenberg group and of their curvatures.

\smallskip

Here is a short overview of the article. 
In Section 2 we provide some background and preliminary results.
In Section 3 we prove Theorem \ref{Tmain}. Section 4 is devoted to a 
discussion of the mean curvature,
while in Section 5 we discuss several examples.

{\bf Note on the proofs.} We will always assume implicitely that our surface $S$ 
is locally parametrizable as $0=g(x,y,t)=t-f(x,y)$. In particular, without loss of generality, we 
assume that $\Omega$ lies below the graph of $f$. The cases when $S$ does not admit a 
parametrization of this kind, i.e., in the ``near-Euclidean case'' 
when the plane tangent to $S$ is vertical at some point,
can be dealt with in a way similar to that employed in \cite{AF1}.

It is a pleasure to thank R. Monti for his constructive criticism and L. Capogna for directing us
to several references.

\section{Notation and preliminaries}\label{definizioni}
In this section, we collect some basic definitions and known facts about the
structure and the geometry of $\Hi$. There is a vast literature on
sub-Riemannian geometry and Carnot groups. Justs to quote a few
titles, we refer the reader to
\cite{B}, \cite{FSSC1}, \cite{G}, \cite{Montgomery}, \cite{NSW},
\cite{Pansu}, \cite{Pansu2}, \cite{Stein}.

The Heisenberg group $\Hi=\Hi^1$ is the Euclidean space ${\mathbb R}^3$
endowed with the noncommutative product
$$
(x,y,t)\cdot(x^\prime,y^\prime,t^\prime)=
(x+x^\prime,y+y^\prime,t+t^\prime+2(x^\prime y-xy^\prime)).
$$
Sometimes it is convenient to think of the elements of  $\Hi$ as $(z,t)\in{\mathbb
  C}\times{\mathbb R}$.
The {\bf Carnot-Charath\'eodory distance} in $\Hi$ is
the sub-Riemannian metric that makes  pointwise
orthonormal the left invariant vector fields
$X$ and $Y$,
$$
X=\partial_x+2y\partial_t,\ Y=\partial_y-2x\partial_t.
$$
The vector fields $X$, $Y$ do not commute, $[X,Y]=-4\partial_t$.
The distance between two points $P$ and $Q$ in $\Hi$ is denoted by $d(P,Q)$.
The span of the vector fields $X$ and $Y$ is called {\bf horizontal
  distribution}, and it is denoted by $\HH$. The
{\bf fiber} of $\HH$ at a point $P$ of $\Hi$ is
$\HH_P=\text{span}\{X_P,Y_P\}$. The inner product in $\HH_P$ is denoted by
$\langle\cdot,\cdot\rangle$.

An important element of $\Hi$'s structure is the {\bf dilation
  group at the origin}
$\{\delta_\lambda: \ \lambda\ne 0\}$,
$$
\delta_\lambda(z,t)=(\lambda z,\lambda^2 t),\ z=x+iy
$$
By left translation, a dilation group is defined at each point $P$ of
$\Hi$.

The Heisenberg group is also endowed with a {\bf rotation group},
which is useful in simplifying some calculations. For
$\theta\in{\mathbb R}$, let
$$
R_\theta(z,t)=(e^{i\theta},t)
$$
be the rotation by $\theta$ around the $t$-axis. Composing with left
translation, one could define rotations around any vertical line
$(x,y)=(a,b)$. $R_\theta$ is an isometry of $\Hi$ and its differential
acts on the fiber ${\mathcal H}_O$ as a rotation by $\theta$. Under the
usual identification between the Riemannian tangent space of $\Hi$ at $O$, $T_O\Hi$,
and the Lie algebra $\mathfrak h$  of $\Hi$,
the differential of $R_\theta$ can be thought of as a rotation on $\mbox{span}\{X,Y\}$, the
first stratum of $\mathfrak h$. With respect to the basis $\{X,Y\}$,
$$
dR_\theta V=
\begin{pmatrix}
\cos(\theta)&-\sin(\theta)\cr

\sin(\theta)&\cos(\theta)
\end{pmatrix}
V
$$
whenever $V\in\mbox{span}\{X,Y\}$.
With some abuse of notation, we denote $dR_\theta$ by $R_\theta$.

The distance between two points $P$ and $Q$ in $\Hi$ is defined as
follows. Consider an absolutely continuous curve $\gamma$ in
${\mathbb R}^3$, joining $P$ and $Q$, which is {\bf horizontal}. That
is, $\dot{\gamma}(t)=a(t)X_{\gamma(t)}+b(t)Y_{\gamma(t)}$ lies in $\HH_{\gamma(t)}$. Its
Carnot-Charath\'eodory length is $\length(\gamma)$,
$$
\length(\gamma)=\int(a(t)^2+b(t)^2)^{1/2}dt
$$
 The
Carnot-Charath\'eodory distance between $P$ and $Q$, $d(P,Q)$,
is the infimum of the Carnot-Charath\'eodory lengths of such
curves. The infimum is actually a minimum, the distance between $P$
and $Q$ is realized by the length of a geodesic. By translation
invariance, all geodesics are left translations of geodesics passing through
the origin.
The unit-speed geodesics at the origin
\cite{Monti}, \cite{Montgomery} are
\begin{equation}
\label{eqGeo}
\gamma_{O,\phi,W}(\sigma)=
\begin{cases}
x(\sigma)=\sin(\alpha(W))\frac{1-\cos(\phi\sigma)}{\phi}+\cos(\alpha(W))\frac{\sin(\phi\sigma)}{\phi}\\
y(\sigma)=\sin(\alpha(W))\frac{\sin(\phi\sigma)}{\phi}-\cos(\alpha(W))\frac{1-\cos(\phi\sigma)}{\phi}\\
t(\sigma)=2\frac{\phi\sigma-\sin(\phi\sigma)}{\phi^2}
\end{cases}
\end{equation}
Here, $W$ is a unitary vector in $\H_O$ and $\alpha(W)\in[0,2\pi)$ is
unique with the property $\dot{\gamma}_{O,\phi,W}(0)=W$.
$\phi\in\RR$, and the geodesic is length
minimizing over any interval of length $2\pi/|\phi|$.
In the case $\phi=0$, the geodesic is a straight line in the plane
$\{t=0\}$,
$$
x(\sigma)=\cos(\alpha(W))\sigma,\ y(\sigma)=\sin(\alpha(W))\sigma,
$$
and we say that the geodesic is {\bf straight}.

From these equations we deduce the parametric equations of the boundary of the ball $B(0,r).$
 $(z,t)\in \partial B(0,r)$ if and only if there is $\phi\in [-2\pi/r,2\pi/r]$ so that
\begin{equation}\label{ballprofile}
\begin{cases}
\vert z\vert=2\frac{\sin(\phi r/2)}{\phi}\\
t=2\frac{\phi r-\sin(\phi r)}{\phi^2}.
\end{cases}
\end{equation}

If $P=(z,t)$ and $z\ne0$, then there exists a unique length
minimizing geodesic connecting $P$ and $O$. If $P=(0,t)$, $t\ne0$,
(i.e., if $P$ belongs to the center of $\Hi$)
then there is a one parameter family of length minimizing geodesics joining $P$
and $O$,  obtained by rotation of a single geodesic around the $t$-axis.

We call {\bf curvature} of $\gamma_{O,\phi,W}$ the parameter
$\kappa(\gamma_{O,\phi,W})=|\phi|$.
We will see later how $\kappa(\gamma_{O,\phi,W})$  is related
to the curvature of a surface in $\Hi$.
The notion of curvature extends to all geodesics by
left translations, $\kappa(L_P\gamma)=\kappa(\gamma)$, where $L_P$ is
left translation by $P\in\Hi$ in $\Hi$.
Given points $P=(z,t)$ and $P^\prime=(z^\prime,t^\prime)$, they are
joined by a unique length minimizing geodesic, unless $z=z^\prime$.

Let
$$
\gamma_{P,\phi,\alpha}=L_P\gamma_{O,\phi,\alpha}
$$
The parameter $\phi$ is geometric in the following
sense. $2\pi/|\phi|$ is the length of $\gamma_{P,\phi,W}$ and
$sgn(\phi)$ is positive if and only if the $t$-coordinate increases with
$\sigma$. Recall that in $\Hi$ the orientation of the $t$-axis is an
intrinsic notion, unlike the Euclidean space.
If $\gamma_{P,\phi,W}$ and
$\gamma_{P^\prime,\phi^\prime,W^\prime}$ have an arc in common,
then $\phi=\pm\phi^\prime$, while no such easy relation exists for the
parameter $W$.
To change the orientation of a geodesic, observe that
$$
\gamma_{P,\phi,W}(\sigma)=\gamma_{P,-\phi,-W}(-\sigma)
$$
Unlike the Euclidean case, a geodesic $\gamma$ leaving $O$ is not
determined by its tangent vector at the origin,
$\dot{\gamma}(0)=W=\cos(\alpha(W))X_O+\sin(\alpha(W))Y_O$. The extra parameter
we need is $\kappa(\gamma)$.

The parameter $\kappa(\gamma)$, if $\phi=\kappa(\gamma)\ne 0$, is related to the dilation group as
follows. $\delta_\lambda\left(\gamma_{O,\phi,W}\right)$ is a
reparametrization of the geodesic $\gamma_{O,\phi/\lambda,W}$. That
is, all geodesics $\gamma$ leaving $0$ and having fixed initial
velocity $\dot{\gamma}(0)=v$ in $\HH_O$ are dilated of each other.
The case of the straight geodesics is the limiting one, corresponding
to $\lambda\to0$. In a precise sense, then, the set of geodesics at
$O$ is parametrized by the unit circle in $\HH_O$ and by the dilation
group, a feature of $\Hi$ with no Euclidean counterpart. This simple
fact will be a recurrent theme in the following sections.

\medskip

Let $S$ be a smooth surface in $\Hi$.
We need some geometric notions about $S$, some of which were studied in \cite{AF1}.
\begin{definition}
Let $E$ be a closed subset of $\Hi$, $P\in E$. The {\bf metric  normal} to
$E$ at $P$
  is the set $\N_PS$ of the points $Q\in\Hi$ such that $d(Q,E)=d(Q,P)$.
\end{definition}
\begin{definition}
Let $S$ be a $C^1$ surface in the Euclidean sense in $\Hi$, 
which is the boundary of an open set
$\Omega$ and of $\Hi-\overline{\Omega}$ and let $P\in S$. The 
{\bf oriented metric normal to $S$ at $P$}, $\N^+_PS$, is the unique parametrization
of $\N_PS$ such that $\delta_S(\N^+_PS(\sigma),P)=\sigma$. 
\end{definition}
\begin{definition}
Let $S$ be a smooth surface in $\Hi$ and
let $P\in S$ be a non-characteristic point.
The Euclidean tangent space to $S$ at $P$ is denoted by $T_PS$.
The {\bf direction tangent to
 $S$ at $P$} is the $1$-dimensional space $V_PS=T_PS\cap \HH_P$.
The {\bf plane tangent to
 $S$ at $P$}, $\Pi_PS$, is the Euclidean plane in $\Hi$, tangent to
$S$ at $P$ in the Euclidean sense. The {\bf direction normal to} $S$
at $P$ is $N_PS=\HH_P\ominus V_PS$, where $\ominus$ is taken with respect to the 
subRiemannian metric.

The {\bf exterior normal to} $S$ at $P$ is 
$$
N_PS=\left({\mathcal N}_PS\right)^\cdot(0),
$$
the vector tangent to ${\mathcal N}_PS$ at $P$.
\end{definition}
The {\bf group tangent to $S$ at $P$} \cite{Pansu}
is the $2$-dimensional vector space $G_PS={\mathcal V}_PS\oplus \T$,
where $\T=\{(0,0,t):\ t\in{\mathbb R}\}$ is the center of $\Hi$ and ${\mathcal V}_PS$ is the one
parameter subgroup of $\Hi$ generated by $V_PS$. One point we want to
make in the present paper is that $G_PS$ does not seem to capture the
complexity of the geodesics' set, while $T_PS$ does, in a precise sense.

The following facts are easily established by direct calculation.
\begin{proposition}
Let $S$ be a smooth surface in $\Hi$, implicitely defined by
$g(x,y,t)=0$. Let $P\in S$ be non-characteristic. Then,
\begin{equation}
\label{eqnortan}
V_PS=\text{span}\{Yg\cdot X-Xg\cdot Y\},\
N_PS=\text{span}\{Xg\cdot X+Yg\cdot Y\}=\text{span}\{\nabla_\Hi g\}
\end{equation}
\end{proposition}
Suppose that $S$ is locally parametrized by $g(x,y,t)=0$. Let $P$ be a point in $S$ and let 
$C$ be the characteristic point of $\Pi_PS$, the euclidean plane in $\Hi$ which is tangent to $S$
at $P$. Then,
\begin{equation}
\label{eqchardist}
d(P,C)=\frac{|\nabla_\Hi g(P)|}{2|\partial_tg(P)|}
\end{equation}

The equation of $\N^+_PS$ can be written in terms of
$g$'s partial derivatives.
$\N^+_PS=
P\cdot\eta$ (left translation by $P$), where
$\eta=(u,v,s)$ and
\begin{equation}
\label{eqnorSo}
\eta(\sigma)=
\begin{cases}
x^\prime(\sigma)=
\frac{1}{4\partial_tg}\left\{Yg(P)
\left(1-\cos\left(\frac{4\partial_tg(P)\sigma}{|\nabla_\Hi g(P)|}\right)\right)
+Xg(P)\sin\left(\frac{4\partial_tg(P)\sigma}{|\nabla_\Hi g(P)|}\right)\right\}\\
y^\prime(\sigma)=
\frac{1}{4\partial_tg}\left\{-Xg(P)
\left(1-\cos\left(\frac{4\partial_tg(P)\sigma}{|\nabla_\Hi g(P)|}\right)\right)
+Yg(P)\sin\left(\frac{4\partial_tg(P)\sigma}{|\nabla_\Hi g(P)|}\right)\right\}\\
t^\prime(\sigma)=\frac{|\nabla_\Hi g(P)|^2}{8(\partial_t g(P))^2}
\left\{\frac{4\partial_tg(P)\sigma}{|\nabla_\Hi g(P)|}
-\sin\left(\frac{4\partial_tg(P)\sigma}{|\nabla_\Hi g(P)|}\right)\right\}
\end{cases}
\end{equation}
See p.10 in \cite{AF1}. When $g(z,t)=t-f(z)$, we have
\ieqas
\left\{\begin{array}{l}
x'=\frac{1}{4}(Xf \sin\alpha + Yf  (1-\cos\alpha))\\
y'=\frac{1}{4}(Yf \sin\alpha - Xf  (1-\cos\alpha))\\
t'=\frac{\mid\nabla_\Hi f\mid^2}{8}(\alpha -\sin \alpha),
\end{array}
\right.,
\eeqas
where $\alpha =\frac{4\tau}{\mid\nabla_{\Hi}f\mid}$.
We introduce an {\bf exponential map} $F:(u,v,\sigma)\mapsto F(u,v,\sigma)$, defined from an open 
subset of ${\mathbb R}^2\times{\mathbb R}$ with values in $\Hi$,
\begin{equation}
\label{eqexp}
F(u,v,\sigma)=\N^+_{(u,v,f(u,v))}S(\sigma).
\end{equation}
We are going to use the coordinates
$$
(x,y,t)=F(u,v,\sigma)=(u,v,f(u,v))\circ (x',y',t')
$$
where $(x',y',t')=\gamma_{u,v}(\sigma)$ and $\gamma_{u,v}$ is the
metric normal's left translate by $P^{-1}$. We have, then
\ieqas
\left\{\begin{array}{l}
x_u=1+x'_u\\
y_u= y'_u\\

t_u=f_u(u,v)-2y' +t'_u+2(vx'_u-uy'_u),
\end{array}
\right.
\eeqas
See \cite{AF1}, p.15, Lemma 6.5.

\smallskip

At p.16 of \cite{AF1} we find the following formulae, which we summarize in a vast lemma.
\begin{lemma}\label{vast} We have the following equalities.
\ieqas
4x'_u& =&(Xf)_u\sin\alpha +(Yf)_u (1-\cos\alpha)+\alpha_u Xf\cos\alpha
+\alpha_u Yf\sin\alpha\\
&=&(Xf)_u\sin\alpha +(Yf)_u
(1-\cos\alpha)\\
&-&\cos\alpha\frac{Xf}{2}\frac{4\tau}{\mid
\nabla_\Hi
f\mid^3}\big(2Xf(Xf)_u+2Yf(Yf)_u\big)-\sin\alpha\frac{Xf}{2}\frac{4\tau}{\mid
\nabla_\Hi f\mid^3}\big(2Xf(Xf)_u+2Yf(Yf)_u\big),
\eeqas
and so $(x'_u)_{\mid \tau =0}=0$.
Analogously
\ieqas
4x'_v& =&(Xf)_v\sin\alpha +(Yf)_v (1-\cos\alpha)\\
&-&\big(\cos\alpha Xf\frac{4\tau}{\mid
\nabla_\Hi f\mid^3}+\sin\alpha Yf \frac{4\tau}{\mid
\nabla_\Hi
f\mid^3})\big ( Xf(Xf)_u+ Yf(Yf)_u\big),
\eeqas
and $(x'_v)_{\mid \tau =0}=0$;
\ieqas
4x'_\tau& =& Xf \cos\alpha \frac{4 }{\mid
\nabla_\Hi f\mid } + Yf  \sin\alpha \frac{4 }{\mid
\nabla_\Hi f\mid },
\eeqas
and
$
(x'_\tau)_{\mid \tau =0}=\frac{Xf}{\mid
\nabla_\Hi f\mid }$;
\ieqas
4y'_u& =&(Yf)_u\sin\alpha -(Xf)_u (1-\cos\alpha)\\
&-& \frac{4\tau}{\mid
\nabla_\Hi
f\mid^3}\big(\cos\alpha Yf-\sin\alpha Xf\big)\big ( Xf(Xf)_u+ Yf(Yf)_u\big),
\eeqas
and
$
(y'_u)_{\mid \tau =0}=0$;
\ieqas
4y'_v& =&(Yf)_v\sin\alpha -(Xf)_v (1-\cos\alpha)\\
&-& \frac{4\tau}{\mid
\nabla_\Hi
f\mid^3}\big(\cos\alpha Yf-\sin\alpha Xf\big)\big ( Xf(Xf)_v+
Yf(Yf)_v\big),
\eeqas
and
$
(y'_v)_{\mid \tau =0}=0$;
\ieqas
4y'_\tau& =&  \frac{4\tau}{\mid
\nabla_\Hi
f\mid^3}\big(\cos\alpha Yf-\sin\alpha Xf\big)  ,
\eeqas
and
$
(y'_\tau)_{\mid \tau =0}=\frac{ Yf}{\mid\nabla_\Hi f\mid }
$.
\ieqas
8t'_u& =&2\big ( Xf(Xf)_u+
Yf(Yf)_u\big)(\alpha-\sin\alpha)-
 \frac{4\tau}{\mid
\nabla_\Hi
f\mid }\big ( Xf(Xf)_u+
Yf(Yf)_u\big)(1-\cos\alpha),
\eeqas
and
$(t'_u)_{\mid \tau =0}=0$;
\ieqas
8t'_v& =&2\big ( Xf(Xf)_v+
Yf(Yf)_v\big)(\alpha-\sin\alpha)-
 \frac{4\tau}{\mid
\nabla_\Hi
f\mid }\big ( Xf(Xf)_v+
Yf(Yf)_v\big)(1-\cos\alpha),
\eeqas
and
$
(t'_v)_{\mid \tau =0}=0,$
and, eventually
\ieqas
8t'_\tau& =&
  4  \mid
\nabla_\Hi
f\mid (1-\cos\alpha),
\eeqas
and
\ieqas
t'_\tau =0.
\eeqas
\end{lemma}
\section{The Hessian of the function measuring the distance from a 
smooth surface}\label{hessian}
Let $g$ be a smooth function from $\Hi$ in ${\mathbb R}$. We define
the {\bf horizontal Hessian} $\mbox{Hess}_\Hi g$ and the {\bf symmetrized horizontal Hessian}
$\mbox{Hess}^*_\Hi g$
of $g$. See, e.g., \cite{DGN}, \cite{LMS} for a discussion of second
order differential operators in the Heisenberg group.
\begin{equation}
\label{eqhessdf}
\mbox{Hess}_\Hi g=
\begin{pmatrix}
XXg & YXg \cr
XYg & YYg
\end{pmatrix}
\end{equation}
Let $(XY)^*g=\frac{1}{2}(XYg+YXg)$.
\begin{equation}
\label{eqhessdfs}
\mbox{Hess}^*_\Hi g=
\begin{pmatrix}
XXg & (XY)^*g \cr
(XY)^*g & YYg
\end{pmatrix}
\end{equation}
Next, we prove Theorem \ref{Tmain}. For ease of calculation, we restate it in a 
less geometric form.
\begin{theorem}
\label{Thess}
Let $S=\partial\Omega$ be a non-characteristic smooth surface in $\Hi$, where $\Omega$ 
locally defined by the equation
$t-f(x,y)<0$.. Let
$$
Q=XXf\ (Yf)^2-2(XY)^*f\ Xf\ Yf+YYf\ (Xf)^2
$$
Then, the horizontal Hessian of ${\delta_S}$ is
\begin{equation}
\label{eqhessns}
\begin{split}
\mbox{Hess}_\Hi\ {\delta_S}&=
\begin{pmatrix}
\frac{(Yf)^2}{|\nabla_\Hi f|^5}Q+4\frac{Xf\ Yf}{|\nabla_\Hi f|^3} &
-Yf\left[ \frac{Xf}{|\nabla_\Hi f|^5}Q-4\frac{Yf}{|\nabla_\Hi f|^3}
\right]\cr
-Xf\left[ \frac{Yf}{|\nabla_\Hi f|^5}Q+4\frac{Xf}{|\nabla_\Hi f|^3} \right]&
\frac{(Xf)^2}{|\nabla_\Hi f|^5}Q-4\frac{Xf\ Yf}{|\nabla_\Hi f|^3}
\end{pmatrix}\\
&\\
&=v_PS\otimes v_PS \cdot (h_S(P)I+p_S(P)J)
\end{split}
\end{equation}
Hence, the symmetrized horizontal Hessian of the distance
function ${\delta_S}$, restricted to the surface $S$, is
\begin{equation}
\label{eqhess}
\mbox{Hess}^*_\Hi\ {\delta_S}=
\begin{pmatrix}
-\frac{(Yf)^2}{|\nabla_\Hi f|^5}Q-4\frac{Xf\ Yf}{|\nabla_\Hi f|^3} &
\frac{Xf\ Yf}{|\nabla_\Hi f|^5}Q+2\frac{(Xf)^2-(Yf)^2}{|\nabla_\Hi f|^3} \cr
\frac{Xf\ Yf}{|\nabla_\Hi f|^5}Q+2\frac{(Xf)^2-(Yf)^2}{|\nabla_\Hi  f|^3}  &
-\frac{(Xf)^2}{|\nabla_\Hi f|^5}Q+4\frac{Xf\ Yf}{|\nabla_\Hi f|^3}
\end{pmatrix}
\end{equation}
and, on $S$,
\begin{equation}
\label{eqcurvQ}
\Delta_\Hi {\delta_S}=\frac{Q}{|\nabla_\Hi f|^3}=
\frac{XXf(Yf)^2-2(XY)^*fXfYf+YYf(Xf)^2}{\mid \nabla_\Hi
f\mid^3}.
\end{equation}
In particular, if $C$ is a characteristic point in $S$, then
$|\mbox{Hess}_\Hi\ {\delta_S}(P)|=O(d(P,C)^{-1})$ as $P\to C$ in $S$.
\end{theorem}
Elementary linear algebra applied to (\ref{eqhessns}), or a direct calculation with (\ref{eqhess}),
 shows that the characteristic polynomial of 
$\mbox{Hess}^*_\Hi\ {\delta_S}$ is
\begin{equation}
\label{convexitas}
{\mathcal P}(\lambda)=\lambda^2-h_S(P)\lambda-\frac{1}{4}p_S^2(P).
\end{equation}
We could consider the roots $\lambda_{1,2}$ of ${\mathcal P}$, i.e. the eigenvalues of the 
symmetrized horizontal Hessian $\mbox{Hess}^*_\Hi\delta_S$, 
to be the ``horizontal principal curvatures'' of $S$ at $P$. It would be interesting to have a 
direct geometric interpretation for these two quantities.

The eigenvalues of $\mbox{Hess}^*_\Hi\delta_S$ have opposite sign, unless $p_S(P)=0$. 
On the other hand, if there is an open subset $A$ of $S$ s.t. $p_S(P)=0$ 
whenever $P\in A$, then $A$ must be 
 a {\bf vertical} subset of $S$, i.e. $A$ is defined as $\{(x,y,t):\ \varphi(x,y)=0\}$.  
Following \cite{DGN} (see also \cite{LMS} and \cite{DGNT}), we say that a function $k$ is
{\bf mildly convex} if $\mbox{Hess}^*_\Hi k$ is positive definite. As a consequence of the above,
we have the following remark.
\begin{remark}
Let $S$ be an orientable surface in $\Hi$,
which is free of charcateristic points, and let $\delta_S$ the signed distance from $S$ 
associated to an orientation of $S$. Then, $\delta_S$ is mildly convex or concave 
in a neighborhood of $S$
if and only if $S$ is a vertical set  $\{(x,y,t):\ \varphi(x,y)=0\}$ and 
$\{(x,y):\ \varphi(x,y)=0\}$ is a curve in ${\mathcal R}^2$ 
which does not change the sign of its Euclidean 
curvature.

In particular the Carnot-Charath\'eodory distance in $\Hi$ {\bf is not} mildly convex, since the
metric spheres are not vertical sets. This could be compared with the result in \cite{DGN}, where
it is proved that the gauge distance in $\Hi$ is mildly convex. 
\end{remark}

\smallskip

In order to prove the theorem, we need some preliminary results.
We denote by ${\mbox{Hess}}\ k$ the Euclidean Hessian of a function $k$
from ${\mathbb R}^3$ in ${\mathbb R}$. If $K=(k_1,\dots,k_m)$ maps
${\mathbb R}^3$ in ${\mathbb R}^m$, we write

$$
{\mbox{Hess}\ K}=({\mbox{Hess}}\ k_1,\dots,{\mbox{Hess}}\ k_m)
$$
Consider the equation
$$
(x,y,t)=F(u,v,\tau)
$$
Let us denote by ${\mbox{Hess}}\ x$,
${\mbox{Hess}}\ y$, ${\mbox{Hess}}\ t$ the Euclidean Hessian matrices of
$x(u,v,\tau)$, $y(u,v,\tau)$, $t(u,v,\tau)$; let  $\xi=(u,v,\tau)$ and
let $\xi_x$, $\xi_y$, $\xi_t$ be the coloumns of the  Jacobian
matrix of $F^{-1}$ computed when $\tau=0$.

Let $a=\xi_x+2y\xi_t$ and $b=\xi_y-2x\xi_t$.
\begin{lemma}
\label{lemmahess} The entries of the horizontal Hessian of ${\delta_S}$ satisy the
following relations.
\begin{enumerate}
\item[(i)] We have
\begin{equation}
\label{eqXX}
|\nabla_\Hi f|\cdot XX {\delta_S}=-\Big\vert F_u\vert F_v\vert {\mbox{Hess}}F(a,a)\Big\vert
\end{equation}
Here, $\Big\vert c_1|c_2|c_3\Big\vert$ denotes the determinant of the
$3\times 3$ matrix having columns $c_j$, $j=1,2,3$.
\item[(ii)] We have
\begin{equation}
\label{eqYY}
|\nabla_\Hi f|\cdot YY {\delta_S}=-\Big\vert F_u\vert F_v\vert {\mbox{Hess}}F(b,b)\Big\vert
\end{equation}
\item[(iii)] $XY\tau$ and $YX\tau$ satisfy
\begin{equation}
\label{eqXY}
\begin{cases}
X{\delta_S}\cdot XX {\delta_S}+Y{\delta_S}\cdot XY{\delta_S}=0\cr
X{\delta_S}\cdot YX{\delta_S}+Y{\delta_S}\cdot YY{\delta_S}=0
\end{cases}
\end{equation}
\end{enumerate}
\end{lemma}

\begin{proof}
The Inverse Function Theorem applied to the equation $(x,y,t)=F(u,v,\tau)$
gives a set of equations for the second
derivatives of $\tau={\delta_S}$, w.r.t. $x,y,t$. Let $p,q$ be any two
variables, possibly equal, chosen among $x,y,t$. Then
\ieqa\label{sist1}
\left\{\begin{array}{l}
x_uu_{pq}+x_vv_{pq}+x_\tau\tau_{pq}=-  {\mbox{Hess}} x (\xi_p,\xi_q)\\
y_uu_{pq}+y_vv_{pq}+y_\tau\tau_{pq}=-  {\mbox{Hess}} y (\xi_p,\xi_q)\\
t_uu_{pq}+t_vv_{pq}+t_\tau\tau_{pq}=-  {\mbox{Hess}} t (\xi_p,\xi_q)
        \end{array}
\right.
\eeqa

We solve the linear systems with respect to the second partial
derivatives of $\tau$,
\ieqa
\label{eqtpq}
\tau_{pq}=\frac{\mbox{det}
\left[\begin{array}{lll}
x_u,&x_v,&- {\mbox{Hess}} x (\xi_p,\xi_q)\\
y_u,&y_v,&- {\mbox{Hess}} y (\xi_p,\xi_q)\\
t_u,&t_v,&- {\mbox{Hess}} t (\xi_p,\xi_q)
  \end{array}
\right]}{\mbox{det}
\left[\begin{array}{lll}
x_u,&x_v,&x_\tau\\
y_u,&y_v,&y_\tau\\
t_u,&t_v,&t_\tau.
  \end{array}
\right]}
\eeqa
From Lemma \ref{vast},
$$
\mbox{det}
\left\vert\begin{array}{lll}
x_u,&x_v,&x_\tau\\
y_u,&y_v,&y_\tau\\
t_u,&t_v,&t_\tau.
  \end{array}
\right\vert_{\tau=0}
=|\nabla_\Hi f|
$$
Since $XX=\partial_{xx}+4y\partial_{xt}+4y^2\partial_{tt}$,  since the
determinant is linear with respect to each column and
${\mbox{Hess}}F(U,V)$ is bilinear in $U$ and $V$, then (\ref{eqtpq})
implies that, when $\tau=0$,
\begin{eqnarray*}
|\nabla_\Hi f|\cdot XX {\delta_S}&
=&\Big\vert F_u\vert F_v\vert F_\tau\Big\vert\cdot XX {\delta_S}\cr
&=&\Big\vert F_u\vert F_v\vert -[{\mbox{Hess}}F(\xi_x,\xi_x)+4y {\mbox{Hess}}F(\xi_x,\xi_t)+
4y^2 {\mbox{Hess}}F(\xi_t,\xi_t)]\Big\vert\cr
&=&-\Big\vert F_u\vert F_v\vert {\mbox{Hess}}F(a,a)\Big\vert
\end{eqnarray*}
and this shows (i). (ii) can be proved in the same way, starting from
the expression
$$
YY=\partial_{yy}-4x\partial_{yt}+4x^2\partial_{tt}.
$$
(iii) follows taking derivatives w.r.t. $X$ and $Y$ of the  eikonal equation (\ref{eikonal}).
\end{proof}

\begin{proof} of  Theorem \ref{Thess}. By (i) and (ii) in Lemma
  \ref{lemmahess}, the expression of $XX{\delta_S}$ and $YY{\delta_S}$ can be
  computed if we know ${\mbox{Hess}}F$ restricted to $\tau=0$ and
$J(F^{-1})(x,y,t)|_{\tau=0}
=(\xi_x\vert \xi_y\vert \xi_t).$ By (iii), this allows us to
  compute $XY{\delta_S}$ and $YX{\delta_S}$ as well.

We apply the classical Implicit Function Theorem to
$$
(x,y,t)=F(u,v,\tau).
$$
Recall from Lemma \ref{vast} that the Jacobian of $F$,
\begin{equation}JF(u,v,0)=J\left(\begin{array}{lll}
x&y&t\\
u&v&\tau
\end{array}\right)_{\tau =0}=
\begin{pmatrix}
1 & 0 & \frac{Xf}{\vert\nabla_\Hi f\vert}\cr
0 & 1 & \frac{Yf}{\vert\nabla_\Hi f\vert}\cr
 f_u& f_v
& \frac{2vXf-2uYf}{\vert\nabla_{\Hi}f\vert}
\end{pmatrix},
\end{equation}
whose determinant is
$$
\mbox{det}JF(u,v,0)=\mid\nabla_\Hi f\mid.
$$
 Its inverse matrix is
\begin{equation}\label{inv}J(F^{-1})(x,y,t)|_{\tau=0}=J\left(\begin{array}{lll}
u&v&\tau\\
x&y&t
\end{array}\right)_{\tau =0}=
\begin{pmatrix}
1+\frac{Xf f_u}{\mid\nabla_\Hi f\mid^2} & \frac{Xf f_v}{\mid\nabla_\Hi
f\mid^2} & -\frac{Xf}{\vert\nabla_\Hi f\vertÒ2}\cr \frac{Yf
f_u}{\mid\nabla_\Hi f\mid^2} & 1+\frac{Yf f_v}{\mid\nabla_\Hi
f\mid^2} &
-\frac{Yf}{\vert\nabla_\Hi f\vert^2}\cr
 -\frac{ f_u}{\mid\nabla_\Hi
f\mid }& -\frac{  f_v}{\mid\nabla_\Hi
f\mid }
&\frac{1}{\mid\nabla_\Hi
f\mid }.
\end{pmatrix}=(\xi_x\vert \xi_y\vert \xi_t)
\end{equation}
If $g(x,y,t)=t-f(x,y)$, we write $XXf=XXg$, $XYf=XYg$, and so on.
The following relations will be used below.
\begin{equation}
\begin{array}{cc}
(Xf)_u=XXf & (Yf)_u=XYf \\
(Xf)_v=YXf & (Yf)_v=YYf
\end{array}
\end{equation}
Recalling the derivatives of $x,y,t$ computed in the proof of Lemma
\ref{vast}, we can compute
\ieqas
x_{uu}&=&\frac 14 (Xf)_{uu}\sin\alpha +\frac 14(Yf)_{uu}
(1-\cos\alpha)+(\frac 14  (Xf)_{u}\cos\alpha +\frac
14(Yf)_{u}\sin \alpha)\alpha_u\\ &+&\frac
14\alpha_u\sin\alpha\frac{Xf}{2}\frac{4\tau}{\mid
\nabla_\Hi
f\mid^3}\big(2Xf(Xf)_u+2Yf(Yf)_u\big)\\
&-&\frac
14\cos\alpha\big (\frac{Xf}{2}\frac{4\tau}{\mid
\nabla_\Hi f\mid^3}\big(2Xf(Xf)_u+2Yf(Yf)_u\big))_u\\
&-&\frac 14 \cos\alpha\frac{Xf}{2}\frac{4\tau}{\mid
\nabla_\Hi
f\mid^3}\big(2Xf(Xf)_u+2Yf(Yf)_u\big)\\
&-&\frac
14\sin\alpha(\frac{Xf}{2}\frac{4\tau}{\mid
\nabla_\Hi f\mid^3}\big(2Xf(Xf)_u+2Yf(Yf)_u\big)\big)_u,
\eeqas
hence
\ieqas
(x_{uu})_{\mid\tau=0}=0.
\eeqas
\ieqas
x_{vu}&=&\frac 14 (Xf)_{vu}\sin\alpha +\frac 14(Yf)_{uv }
(1-\cos\alpha)+( \frac 14 (Xf)_{u}\cos\alpha +\frac 14(Yf)_{u}
\sin\alpha)\alpha_v\\
&+&\frac 14\alpha_v\sin\alpha\frac{Xf}{2}\frac{4\tau}{\mid
\nabla_\Hi
f\mid^3}\big(2Xf(Xf)_u+2Yf(Yf)_u\big)\\
&-&\cos\alpha \frac 14(\frac{Xf}{2}\frac{4\tau}{\mid
\nabla_\Hi
f\mid^3}\big(2Xf(Xf)_u+2Yf(Yf)_u\big))_v\\
&&
-\alpha_v\frac
14\cos\alpha\frac{Xf}{2}\frac{4\tau}{\mid
\nabla_\Hi f\mid^3}\big(2Xf(Xf)_u+2Yf(Yf)_u\big)\\
&-&\frac
14\sin\alpha(\frac{Xf}{2}\frac{4\tau}{\mid
\nabla_\Hi f\mid^3}\big(2Xf(Xf)_u+2Yf(Yf)_u\big))_v
\eeqas
\ieqas
(x_{vu})_{\mid\tau=0}=0.
\eeqas
\ieqas
x_{vv}& =&\frac 14 (Xf)_{vv}\sin\alpha +\frac 14 (Yf)_{vv}
(1-\cos\alpha)\\
 &+&\frac 14 \big( (Xf)_v\cos\alpha +\frac 14 (Yf)_v
\sin\alpha\big )\\
&-&(\frac 14 \big(\cos\alpha Xf\frac{4\tau}{\mid
\nabla_\Hi f\mid^3}
+\frac 14\sin\alpha Yf \frac{4\tau}{\mid
\nabla_\Hi
f\mid^3})\big ( Xf(Xf)_u+ Yf(Yf)_u\big))_v,
\eeqas
\ieqas
(x_{vv})_{\mid\tau=0}=0.
\eeqas
\ieqas
 x_{u\tau}& =& \frac 14 (Xf\frac{4 }{\mid
\nabla_\Hi f\mid })_u \cos\alpha-\frac 14 Xf\frac{4 }{\mid
\nabla_\Hi f\mid } \sin\alpha \alpha_u\\
&+&
\frac 14( Yf  \frac{4 }{\mid
\nabla_\Hi f\mid }  )_u\sin\alpha
 +  \frac 14 Yf\alpha_u  \cos\alpha \frac{4 }{\mid
\nabla_\Hi f\mid },
\eeqas
\ieqas
(x_{u\tau})_{\mid\tau=0}=\frac{Yf}{\mid \nabla_\Hi
f\mid^3}((Xf)_uYf-Xf(Yf)_u)=\frac{Yf}{\mid \nabla_\Hi
f\mid^3}((XXf)Yf-Xf(XYf)).
\eeqas
\ieqas
x_{\tau\tau}& =& (Xf\frac{1 }{\mid
\nabla_\Hi f\mid } )_\tau \cos\alpha
-Xf \alpha_\tau\sin\alpha \frac{1 }{\mid
\nabla_\Hi f\mid }\\
&+&(Yf\frac{1 }{\mid
\nabla_\Hi f\mid })_\tau  \sin\alpha
+ Yf \alpha_\tau \cos\alpha \frac{1 }{\mid
\nabla_\Hi f\mid }\\
&=&-Xf \alpha_\tau\sin\alpha \frac{1 }{\mid
\nabla_\Hi f\mid }+Yf \alpha_\tau \cos\alpha \frac{1 }{\mid
\nabla_\Hi f\mid },
\eeqas
and
\ieqas
(x_{\tau\tau})_{\mid\tau=0}=\frac{4Yf}{\mid\nabla_\Hi f\mid^2}.
\eeqas

\ieqas
x_{v\tau}& =& (Xf\frac{1 }{\mid
\nabla_\Hi f\mid } )_v \cos\alpha
-Xf \alpha_v\sin\alpha \frac{1 }{\mid
\nabla_\Hi f\mid }\\
&+&(Yf\frac{1 }{\mid
\nabla_\Hi f\mid })_v  \sin\alpha
+ Yf \alpha_v \cos\alpha \frac{1 }{\mid
\nabla_\Hi f\mid }\\
&=&(Xf\frac{1 }{\mid
\nabla_\Hi f\mid } )_v \cos\alpha +(Yf\frac{1 }{\mid
\nabla_\Hi f\mid })_v  \sin\alpha ,
\eeqas
and
\ieqas
(x_{v\tau})_{\mid\tau=0}=\frac{Yf}{\mid \nabla_\Hi
f\mid^3}((Xf)_vYf-Xf(Yf)_v)=\frac{Yf}{\mid \nabla_\Hi
f\mid^3}((YXf)Yf-Xf(YYf)).
\eeqas
Moreover
\ieqas
 t_{uu}& =& f_{uu}+(\frac 14\big ( Xf(Xf)_u+
Yf(Yf)_u\big)_u(\alpha-\sin\alpha)\\
&+&
(\frac 14\big ( Xf(Xf)_u-
Yf(Yf)_u\big)^2(1-\cos\alpha)\frac{\tau}{|\nabla_\Hi f|^3}\\
&-& \frac{\tau}{2\mid
\nabla_\Hi
f\mid }\big ( Xf(Xf)_u\\
&+&
Yf(Yf)_u\big)(1-\cos\alpha))_u
+\left(\frac{\tau}{2\mid
\nabla_\Hi
f\mid }\big ( Xf(Xf)_u+
Yf(Yf)_u\big)\right)_u\cdot(1-\cos\alpha))\\
&+&2vx'_{uu}-2uy'_{uu}-2y'_{u}-2y'_u,
\eeqas
so that
\ieqas
(t_{uu})_{\tau=0}=(f_{uu}+2vx'_{uu}-2uy'_{uu}-2y'_{u}-2y'_u)_{\tau=0}=f_{uu};
\eeqas
analogously
\ieqas
(t_{vu})_{\tau=0}=f_{vu},
\eeqas
and
\ieqas
(t_{vv})_{\tau=0}=f_{vv},
\eeqas
and
\ieqas
(t_{u\tau})_{\tau=0}=2vx^\prime_{u\tau}-2uy^\prime_{u\tau}-2y^\prime_{\tau},
\eeqas
and
\ieqas
(t_{v\tau})_{\tau=0}=2vx^\prime_{v\tau}-2uy^\prime_{v\tau}+2x^\prime_{\tau},
\eeqas
and
\ieqas
(t_{\tau\tau})_{\tau=0}=2vx'_{\tau\tau}-2uy'_{\tau\tau}.
\eeqas

Hence
\begin{equation*}{\mbox{Hess}} x_{\tau =0}=\frac{Yf}{\mid \nabla_\Hi
f\mid^3}
\begin{pmatrix}
0 & 0 &(XXf)Yf-Xf(XYf) \cr 0 & 0&
(YXf)Yf-Xf(YYf)\cr
 (XXf)Yf-Xf(XYf)&(YXf)Yf-Xf(YYf)
&4\mid \nabla_\Hi
f\mid
\end{pmatrix}.
\end{equation*}
Arguing in the same way we get
\begin{equation*}{\mbox{Hess}} y_{\tau =0}=\frac{Xf}{\mid \nabla_\Hi
f\mid^3}
\begin{pmatrix}
0 & 0 &(XYf)Xf-Yf(XXf) \cr 0 & 0&
(YYf)Xf-Yf(YXf)\cr
 (XYf)Xf-Yf(XXf)&(YYf)Xf-Yf(YXf)
&-4\mid \nabla_\Hi
f\mid
\end{pmatrix}.
\end{equation*}
\begin{equation*}{\mbox{Hess}} t_{\tau =0}=
\begin{pmatrix}
f_{uu} & f_{vu} &2vx'_{u\tau}-2uy'_{u\tau}-2y'_\tau\cr
 f_{vu} & f_{vv}& 2vx'_{v\tau}-2uy'_{v\tau}+2x'_{\tau}\cr
2vx'_{u\tau}-2uy'_{u\tau}-2y'_\tau& 2vx'_{v\tau}-2uy'_{v\tau}+2x'_{\tau}&2vx'_{\tau\tau}-2uy'_{\tau\tau}
\end{pmatrix}.
\end{equation*}
The horizontal Hessian matrix can be written now directly.

A direct calculation shows that
\ieqa
a=\left(\frac{(Yf)^2}{|\nabla_\Hi f|^2},
-\frac{Xf\cdot Yf}{|\nabla_\Hi f|^2},\frac{Xf}{|\nabla_\Hi f|}\right)
\eeqa
(viewed as a column vector) and
\ieqa
b=\left(-\frac{Xf\cdot Yf}{|\nabla_\Hi f|^2},\frac{(Xf)^2}{|\nabla_\Hi    f|^2},
\frac{Yf}{|\nabla_\Hi f|}\right)
\eeqa
We then compute
$$
{\mbox{Hess}}|_{\tau=0}x(a,a)=2\frac{(Yf)^2Xf}{|\nabla_\Hi f|^6}Q+4\frac{(Xf)^2Yf}{|\nabla_\Hi f|^4}
$$
$$
{\mbox{Hess}}|_{\tau=0}y(a,a)=-2\frac{(Xf)^2Yf}{|\nabla_\Hi f|^6}Q-4\frac{(Xf)^3}{|\nabla_\Hi f|^4}
$$
$$
{\mbox{Hess}}|_{\tau=0}t(a,a)=-2\frac{(Yf)^2}{|\nabla_\Hi  f|^4}Q
-4\frac{Xf\cdot Yf}{|\nabla_\Hi f|^2}
+8\frac{(Xf)^2}{|\nabla_\Hi  f|^4}(vYf+uXf)
+4\frac{Xf\cdot Yf}{|\nabla_\Hi f|^6}(vYf+uXf)Q
$$
Recall now that
\begin{eqnarray*}
-|\nabla_\Hi f|\cdot XX {\delta_S}&=&
\left|
\begin{array}{ccc}
1 & 0 & {\mbox{Hess}}x|_{\tau=0}(a,a) \\
0 & 1 & {\mbox{Hess}}y|_{\tau=0}(a,a) \\
2v-Xf & -2u-Yf & {\mbox{Hess}}t|_{\tau=0}(a,a)
\end{array}
\right|\\
&=&-\frac{(Yf)^2}{|\nabla_\Hi f|^4}Q-4\frac{Xf\cdot Yf}{|\nabla_\Hi f|^2}
\end{eqnarray*}
Analogous calculations yield the value of $YY{\delta_S}$.

In order to compute the mixed derivatives, we use (\ref{eqXY}).

Finally, observe that, if $S$ is a $C^3$ surface, then
$$
|{\mbox{Hess}}{\delta_S}|_{\tau=0}=O(|\nabla_\Hi f|^{-1})=O(d(P,C)^{-1})
$$
as $P$ tends to a characteristic point $C$ of $S$.
Let us recal that
$$
Q=XXf\ (Yf)^2-2(XY)^*f\ Xf\ Yf+YYf\ (Xf)^2
$$
Then, the horizontal Hessian of ${\delta_S}$ is that given in (\ref{eqhessns}),
$$
-\mbox{Hess}_\Hi\ {\delta_S}=
\begin{pmatrix}
-\frac{(Yf)^2}{|\nabla_\Hi f|^5}Q-4\frac{Xf\ Yf}{|\nabla_\Hi f|^3} &
Yf\left[ \frac{Xf}{|\nabla_\Hi f|^5}Q-4\frac{Yf}{|\nabla_\Hi f|^3}
\right]\cr
Xf\left[ \frac{Yf}{|\nabla_\Hi f|^5}Q+4\frac{Xf}{|\nabla_\Hi f|^3} \right]&
-\frac{(Xf)^2}{|\nabla_\Hi f|^5}Q+4\frac{Xf\ Yf}{|\nabla_\Hi f|^3}
\end{pmatrix}
$$
i.e.
\begin{equation}\begin{split}
-\mbox{Hess}_\Hi\ {\delta_S}=\frac{Q}{\mid\nabla_\Hi f\mid^3}\left(\begin{array}{cc}
-\frac{(Yf)^2}{\mid\nabla_\Hi f\mid^2} &
 \frac{Xf Yf}{|\nabla_\Hi f|^2}\\ 
  \frac{Xf Yf}{|\nabla_\Hi f|^2}&
-\frac{(Xf)^2}{|\nabla_\Hi f|^2}
\end{array}\right)+\frac{4}{\mid\nabla_\Hi f\mid}\left(\begin{array}{cc}
-\frac{Xf Yf}{\mid\nabla_\Hi f\mid^2} &
 -\frac{(Yf)^2}{|\nabla_\Hi f|^2}\\ 
  \frac{(Xf)^2}{|\nabla_\Hi f|^2}&
 \frac{Xf Yf}{|\nabla_\Hi f|^2}.
\end{array}\right)
\end{split}
\end{equation}
Now, let
\begin{equation*} 
M=\left(\begin{array}{cc}
\frac{(Yf)^2}{\mid\nabla_\Hi f\mid^2} &
 -\frac{Xf Yf}{|\nabla_\Hi f|^2}\\ 
  -\frac{Xf Yf}{|\nabla_\Hi f|^2}&
\frac{(Xf)^2}{|\nabla_\Hi f|^2}
\end{array}\right).
\end{equation*}
Then, if we set $p=\frac{4}{\mid\nabla_\Hi f\mid}$ and 
\begin{equation*}
J=\left(\begin{array}{cc}0&1\\
-1&0\end{array}\right),
\end{equation*}
we have
\begin{equation}\begin{split}
\mbox{Hess}_\Hi\ {\delta_S}=kM+pMJ=M(kI+pJ)=
M\left(\begin{array}{cc}k&p\\
-p&k\end{array}\right).
\end{split}
\end{equation}
Let us remark that keeping in mind (\ref{eqchardist}) and that $g=t-f(x,y)$ we have that
$$
d(P,C)=\frac{  |\nabla_\Hi f| }{2},
$$
so that
$$
p=\frac{2}{d(P,C)},
$$
that represents indeed a curvature.
Moreover let us define 
$$
\cos\theta =\frac{Xf}{|\nabla_\Hi f|},
$$
and
$$
\sin\theta =\frac{Yf}{|\nabla_\Hi f|}.
$$
Then
\begin{equation*} 
M=\left(\begin{array}{cc}
\sin^2\theta &
 -\sin\theta \cos\theta\\ 
  -\sin\theta \cos\theta&
\cos^2\theta
\end{array}\right)=v_PS\otimes v_PS,
\end{equation*}

\end{proof}

The estimate for the Hessian of ${\delta_S}$ on $S$ cannot be extended
outside $S$. In Section \ref{examples} , we see that, for instance, when
$S$ is the plane $t=0$, the behavior of $\mbox{Hess}_\Hi\ {\delta_S}$ is worst
when we approach $0$ near the line $x=y=0$.

\section{Mean Curvature}\label{secmean}
In this section, we make some considerations about the notion of mean curvature for
a surface in $\Hi$ (see, e.g.. \cite{GP}).
\begin{definition}
Let $S=\partial\Omega$ be a smooth surface in $\Hi$, with $\Omega$ 
locally defined by the inequality $g(x,y,t)<0$,.
 The {\bf mean curvature} of $S$ at $P\in S$ is
\begin{equation}
h_S(P)=X\left(\frac{Xg(P)}{|\nabla_\Hi g(P)|}\right)+Y\left(\frac{Yg(P)}{|\nabla_\Hi g(P)|}\right)
\end{equation}
\end{definition}
The curvature is well defined for non characteristic points only.

Here are some characterizations of mean curvature. The first is in
terms of the distance function and it is an immediate consequence of
Theorem \ref{Thess}.
\begin{corollary}
\label{Tcurv}
Let $S=\partial\Omega$ be the boundary of an open, $C^3$ set in
$\Hi$. Then, for any non characteristic point $P$ in $S$,
\begin{equation}
\Delta_\Hi {\delta_S}(P)=h_S(P)
\end{equation}
\end{corollary}
R. Monti pointed out to us that the same result can be obtained by an
explicit calculation, which does not rely on the knowledge of the distance's
Hessian. In \cite{GP}, it is shown that the mean curvature is the trace of a certain
matrix ${\mathcal M}$, which is not the Hessian of the distance function.

The second is in terms of horizontal curves on $S$. First, we state an
immediate consequence of the existence theorem for O.D.E.'s.
\begin{lemma}
\label{oggi}
Let $S$ be a smooth surface in $\Hi$ and let $P\in S$ be
non-characteristic. Then, modulo restrictions and reparametrizations, there exists a unique
horizontal curve
$\Gamma$ through $P$ on $S$. Suppose that $S$ is defined by
$g(x,y,t)=0$.
A parametrization of the curve is the solution of the O.D.E.'s
system
\begin{equation}
\label{eqleaf}
\begin{cases}
\dot{\Gamma}(\tau)=Yg(\Gamma(\tau))\cdot X_{(\Gamma(\tau))}-Xg(\Gamma(\tau))\cdot Y_{(\Gamma(\tau))}\\
\Gamma(0)=P
\end{cases}
\end{equation}
\end{lemma}
Before we state the theorem relating the horizontal curves $\Gamma$ and the mean curvature
of $S$, we fix some notation. We give ${\mathcal H}=span\{X,Y\}$ the orientation for which the oriented angle from $X$ to $Y$ has amplitude $\pi/2$. Let $P\in S=\partial\Omega$. The {\bf positive
horizontal direction $V_PS$} tangent to $S$ is $V_P^+S$, the unit vector in the direction $V_PS$
such that the angle from $N_PS$ to $V_P^+S$ has amplitude $-\pi/2$. i.e. 
$$
V_P^+S=dR_{-\pi/2}N_PS,
$$
$dR_{-\pi/2}$ denoting the differential of $R_{-\pi/2}$. If the curve $\Gamma$ of Lemma \ref{oggi}
is oriented in such a way that $\dot{\Gamma}$ has the same orientation of $V_P^+S$, we say that 
$\Gamma$ is a {\bf positively oriented horizontal curve on $S$}. 
Let $p$ be
the projection of $\Hi$ onto the $z$-plane, $p(z,t)=z$ and let $\gamma=p(\Gamma)$ be the 
projection of $\Gamma$. $\gamma$ inherits the orientation of $\Gamma$. The {\it Euclidean curvature
$k$ of the oriented curve} $\gamma$ at $\gamma(t_0)$ is defined (locally in $t_0$) by
$$
\frac{d}{dt}\vert_{t=t_0}\left(n(\gamma(t))\right)
=k\cdot \dot{\gamma}(t_0),
$$
where
$$
n(\gamma(t))=\frac{dR_{\pi/2}\dot{\gamma}(t)}{|dR_{\pi/2}\dot{\gamma}(t)|}
$$
is the {\it unit normal} to the oriented curve $\gamma$ at $\gamma(t)$, which is locally
well defined even if $\gamma$ has self intersections.

In a picture, $\gamma$ has positive curvature if it borders a convex region on its right hand side.

\begin{theorem}
\label{Tcurveuc}
Let $S=\partial\Omega$ be a smooth surface and let $P\in S$ be non-characteristic. Let
$\Gamma$ be the oriented horizontal curve on $S$ through $P$.  
Then, $h_S(P)$ is the curvature
of the oriented curve $p(\Gamma)$ in $\RR^2$ at $p(P)$.
\end{theorem}
\begin{proof}
Suppose first that $T_PS$ is a non-vertical plane, that is, that in a neighborhood of $P$, $S$ 
admits the
parametrization $S=\{f(x,y)=t\}$. The horizontal curve $\Gamma$ on
$S$ through $P$ satisfies equation (\ref{eqleaf}), hence
$\Gamma=(x,y,t)$ is a solution of
$$
\begin{cases}
\dot{x}=Yf(\Gamma)\cr
\dot{y}=-Xf(\Gamma)\cr
\dot{t}=2yYf(\Gamma)+2xXf(\Gamma)
\end{cases}
$$
Since $Xf$ and $Yf$ are independent of $t$, we solve for $x$ and $y$,
independently of $t$, and we find the differential equation for
$\gamma=p(\Gamma)$,
$$
\begin{cases}
\dot{x}=Yf(\Gamma)\cr
\dot{y}=-Xf(\Gamma)
\end{cases}
$$
Suppose that $Yf\ne 0$ at $P$ (otherwise, we can let $Xf\ne 0$ at $P$,
since $P$ is non characteristic).
Parametrizing $\gamma$ as $y=\gamma(x)$, we find
$$
y^\prime=-\frac{Xf}{Yf}
$$
The curvature of $\gamma$ in the Euclidean plane is defined as
\begin{eqnarray*}
k&=&-\frac{d}{dx}\left(\frac{y^\prime}{\sqrt{1+(y^\prime)^2}}\right)=
-\frac{y^{\prime\prime}}{(1+(y^\prime)^2)^{3/2}}\cr
&=&\frac{(Yf)^3}{|\nabla_\Hi
  f|^3}\frac{d}{dx}\left(\frac{Xf}{Yf}\right)\cr
&=&\frac{Yf}{|\nabla_\Hi  f|^3}
\left[(\partial_xXf+\partial_yXfy^\prime)Yf-(\partial_xYf+\partial_yYfy^\prime)Xf\right]\cr
&=&\frac{XXf\ (Yf)^2-2(XY)^*f\ Xf\ Yf+YYf\ (Xf)^2}{|\nabla_\Hi
  f|^3}\cr
&=&\Delta_\Hi {\delta_S}=h_S(P)
\end{eqnarray*}
by Corollary \ref{Tcurv}.

The case when $T_PS$ is a vertical plane can be dealt with similarly.

\end{proof}

\smallskip

In the Heisenberg group, vectors belonging to different sections of
the tangent bundle can be compared by means of left
translations. Under the identification $\Hi\equiv{\mathbb R}^3$, hence
of the tangent space at $P\in\Hi$ with ${\mathbb R}^3$, the differential
of the left translation by $Q=(a,b,c)$, ${\mathcal L}_Q:P\mapsto
Q\circ P$, is the matrix
$$
d{\mathcal L}_Q=
\left[\begin{array}{ccc}
                        1&0&0\\
0&1&0\\
                        2b&-2a&1
                        \end{array}\right]
$$
which is independent of the point where the differential is calculated.
By definition of Carnot-Carath\'eodory distance, $d{\mathcal L}_Q$
maps ${\mathcal H}_P$ in ${\mathcal H}_{Q\circ P}$, isometrically.

Fix $P\in\Hi$ and let $\phi:I\to\Hi$ be a smooth, horizontal curve,
such that $\phi(0)=P$. Let $W$ be a section of the horizontal bundle
along $\phi$: $W(s)\in{\mathcal H}_{\phi(s)}$ for $s\in I$, and $W$ is
a smooth map.
We define the {\bf derivative} of $W$ along
$\phi$ at $P$ as
$$
\frac{d_\phi^\Hi W(s)}{ds}_
{\vert s=0}:=\lim_{s\to 0}\frac{d{\mathcal L}_{ (
\phi(0))\circ (- \phi(s))} W(s)-W(0)}{s}
$$
 We give the notion of derivative of a vector field
with respect to a vector in the Heisenberg group.
\begin{definition}
Let $W$ be a section of the horizontal bundle ${\mathcal H}$. Let
$V\in{\mathcal H}_P$ and $\phi:I\to \Hi$ be any horizontal curve such
that $\phi(0)=P$ and $\phi^\prime (0)=V$. We set
$$
\nabla^\Hi_V W(P):=\frac{d^\Hi (W\circ\phi)(s)}{ds}_{\vert
s=0}.
$$
\end{definition}
\begin{lemma}
\label{lemmacurvinv}
Let $W=aX+bY$ be a section of the horizontal bundle ${\mathcal H}$
and $\phi:I\to \Hi$ be any horizontal curve such
that $\phi(0)=P$ and $\phi^\prime (0)=V$. Then
$$
\nabla^\Hi_V W(P)=Va(P)X(P)+Vb(P)Y(P)
$$
where $Va$ and $Vb$ are the derivatives of the functions $a$ and $b$
in the direction $V$, respectively.
\end{lemma}
\begin{proof}
Indeed
\ieqas
&&
\frac{d{\mathcal L}_{ (
\phi(0))\circ (- \phi(s))} W(s)-W(0)}{s}\\
&=&\frac{d{\mathcal L}_{ ( \phi(0))\circ (-
\phi(s))}\big(a(\phi(s)),b(\phi(s),2\phi_2(s)a(\phi(s))-2\phi_1(s)b(\phi(s))\big)
-W(\phi(0))}{s}\\
&=&\frac{
\big(a(\phi(s)),b(\phi(s),2\phi_2(0)a(\phi(s))-2\phi_1(0)b(\phi(s))\big)
-W(\phi(0))}{s}\\
&=&\frac{a(\phi(s))-a(\phi(0))}{s}X(\phi(0))+\frac{b(\phi(s))-b(\phi(0))}{s}Y(\phi(0))
 \eeqas
and as a consequence
$$
\lim_{s\to 0}\frac{d{\mathcal L}_{ (
\phi(0))\circ (- \phi(s))} W(\phi(s))-W(\phi(0))}{s}=Va(P)X(P)+Vb(P)Y(P)
$$
\end{proof}

Let $N_S(P)$ the unit vector normal to $S=\partial\Omega$ at $P$, in the
Heisenberg sense, and
pointing outside $\Omega$. If $g(x,y,t)=0$ is a local parametrization
of $S$, then
$$
N_S(P)=
\left(\frac{Xg(P)}{|\nabla_\Hi g(P)|}X_P+\frac{Yg(P)}{|\nabla_\Hi g(P)|}Y_P\right)
$$
depending on the orientation of $S$.
\begin{theorem}
\label{Tcurvinv}
The linear operator $M_PS:V\mapsto \nabla^\Hi_{V}N_S(P)$ maps $V_PS$ into
$V_PS$. $M_PS$ acts on $V_PS$ as multiplication times $h_S(P)$.
\end{theorem}
Modulo a sign, the map $M_PS$ is the Weingarten map for surfaces in the Heisenberg group.
\begin{proof}
If $W=N_S(P)$, we have
$$
\nabla^\Hi_{V}N_S(P)\in V_PS
$$
In  fact, differentiating $|N_S(P)\circ\phi(s)|^2=1$ with respect to
$s$, $N_S(P)=\nu_1(P)X_P+\nu_2(P)Y_P$,
we obtain
$$
0=\nu_1(\phi(s))\frac{d}{ds}\nu_1(\phi(s))+\nu_2(\phi(s))\frac{d}{ds}\nu_2(\phi(s))
$$
For $s=0$, this relation becomes
$$
0=\langle N_S(P),\nabla^\Hi_{V}N_S(P)\rangle_\Hi
$$
Since
$N_S(P)=\frac{Xg(P)}{\vert\nabla_\Hi g(P)\vert}X_P+
\frac{Yg(P)}{\vert\nabla_\Hi g(P)\vert}Y_P$ we get, for
$V=\alpha_1X(P)+\alpha_2Y(P)\in V_PS$,
\begin{equation}
\label{eqdozza}
\nabla^\Hi_{V}
N_S(P)=\left[\begin{array}{cc}X(\frac{Xg(P)}{\mid
\nabla_\Hi g(P)\mid})& Y(\frac{Xg(P)}{\mid
\nabla_\Hi g(P)\mid})\\
X(\frac{Yg(P)}{\mid
\nabla_\Hi g(P)\mid})& Y(\frac{Yg(P)}{\mid
\nabla_\Hi g(P)\mid})
     \end{array}\right]\cdot
\left[\begin{array}{c}\alpha_1\\
 \alpha_2\end{array}\right]
\end{equation}

The linear map $V\mapsto\nabla^\Hi_{V}N_S(P)$ is defined on the
one-dimensional linear space $V_PS$, hence it acts as
multiplication by a constant $k$. Testing (\ref{eqdozza}) on
$(\alpha_1,\alpha_2)=(-Yg(P),Xg(P))=\dot{\Gamma}(t)$, if $\Gamma(t)=P$, 
i.e. on a basis of $V_PS$, we see that
$$
k=\frac{Q}{|\nabla_\Hi g|^3}
$$
is the mean curvature of $S$ at $P$.
\end{proof}

Observe that
multiplication by the matrix in (\ref{eqdozza}) acts on $V_PS$ as multiplication by the matrix' 
trace.

\section{Examples}

\label{examples}
In this section, we see two basic examples of surfaces in
$\Hi$. First, we consider the case of the plane $\{t=0\}$. Then, we
consider the case of the ball with respect to the Heisenberg
distance. In particular, we explicitely compute the Hessian of the
function "distance from a point''.

Calculations are easier if we exploit the symmetries of $\Hi$.
\begin{lemma}
\label{lemmahomog}
Let $g$ be a smooth function on $\Hi$, which is homogeneous of degree
$m$ with respect to the dilation group at
the origin. Let $W_j$ be left invariant vector fields on $\Hi$,
$j=1,\dots,n$. Then, $W_1\dots W_ng$ is homogeneous of degree $m-n$,
\begin{equation}
W_1\dots W_ng(\delta_\lambda P)=\lambda^{m-n}g(P)
\end{equation}
In particular, if $S$ is a smooth surface in $\Hi$ and $\sigma$
denotes the distance function from $S$, then $\sigma$, $\nabla_\Hi\sigma$
and $\mbox{Hess}_\Hi\ \sigma$ are homogeneous of degreee $1$, $0$ and $-1$, respectively.
\end{lemma}
Recall that $R_\theta$ denotes both the rotation by $\theta$ around
the $t$-axis in $\Hi$ and the rotation by $\theta$ in
$\mbox{span}\{X,Y\}$. We use the same symbol for the matrix
$$
R_\theta=
\begin{pmatrix}
\cos(\theta)&-\sin(\theta)\cr
\sin(\theta)& \cos(\theta)
\end{pmatrix}
$$
\begin{lemma}
\label{lemmarot}
Let $f$ be a smooth function on $\Hi$, which is invariant under the
rotation group at the origin,
$$
f(r_\theta P)=f(P)
$$
for $\theta$ in ${\mathbb R}$. Then, $\nabla_\Hi$ and $\mbox{Hess}_\Hi$
commute with rotations.  Namely,
\begin{equation}
\label{eqgradrot}
\nabla_\Hi f(R_\theta P)=R_\theta\nabla_\Hi f(P)
\end{equation}
and
\begin{equation}
\label{eqhessrot}
\mbox{Hess}_\Hi f(R_\theta P)=R_\theta\mbox{Hess}_\Hi f(P)
\end{equation}
\end{lemma}
In (\ref{eqhessrot}), the multiplication on the right hand side is a
multiplication of $2\times2$ matrices.

\subsection{Nonvertical planes}
Let $\Pi$ be the plane $\{t=0\}$. Observe that all nonvertical planes
are isometric to $\Pi$. In \cite{AF1} we computed the
metric normal to $\Pi$ at $P=(x,y,0)=(z,0)$,
$$
\gamma_P(\sigma)=(u(\sigma),\ v(\sigma),\ s(\sigma))=
\begin{cases}
\frac{x}{2}\left(1+\cos\left(\frac{2\sigma}{|z|}\right)\right)+\frac{y}{2}\sin\left(\frac{2\sigma}{|z|}\right)\\
\frac{y}{2}\left(1+\cos\left(\frac{2\sigma}{|z|}\right)\right)-\frac{x}{2}\sin\left(\frac{2\sigma}{|z|}\right)\\
\frac{|z|^2}{2}\left(\frac{2\sigma}{|z|}+\sin\left(\frac{2\sigma}{|z|}\right)\right)
\end{cases}
$$
The geodesic $\gamma_P$ realizes the distance from $\Pi$ when
$|\sigma|\le\frac{\pi|z|}{2}$.
By Lemma \ref{lemmahomog}, $\delta_\Pi$ is homogeneous of degree $1$,
$\nabla_\Hi \delta_\Pi$ is homogeneous of degree $0$ and $\mbox{Hess}_\Hi \delta_\Pi$
is homogeneous of degree $-1$. Since $\delta_\Pi$ is invariant under
rotations around the $t$-axis, $\mbox{Hess}_\Hi \delta_\Pi$
commutes with the same group of rotations.

Let $\beta=\frac{\sigma}{|z|}\in(-\pi/2,\pi/2)$, $R=|z|\ge0$ and $\alpha\in[0,2\pi)$.
Here, we use on $\Hi$ the coordinates
$$
(u,v,s)=(R\cos\beta\cos\alpha, R\cos\beta\sin\alpha,
R^2(\beta+1/2\sin(2\beta))).
$$
Hence
$$
J\left(\begin{array}{ccc}
u,&v,&s\\

R,&\alpha,&\beta
\end{array}\right)=\left[\begin{array}
{ccc}
\cos\alpha\cos\beta,&-R\cos\beta \sin \alpha,&-R\sin\beta\cos\alpha\\
\cos\beta\sin\alpha,&R\cos\beta\cos\alpha,&-R\sin\beta \sin \alpha\\
2R(\beta +\frac{1}{2}\sin (2\beta)),&0,&2R^2\cos^2\beta
                   \end{array}\right]
$$
and $$\mbox{det}J\left(\begin{array}{ccc}
u,&v,&s\\
R,&\alpha,&\beta
\end{array}\right)=2R^3\cos\beta (\cos\beta +\beta \sin\beta).$$

As a consequence
$$\mbox{det}J\left(\begin{array}{ccc}
R,&\alpha,&\beta\\
u,&v,&s
\end{array}\right)=\left[\begin{array}
{ccc}
\frac{\cos^2\beta\cos \alpha}{\cos\beta+\beta
\sin\beta},&\frac{\cos^2\beta\sin \alpha}{\cos\beta+\beta
\sin\beta},&\frac{\sin\beta}{2R(\cos\beta+\beta \sin\beta)}\\
-\frac{\sin\alpha}{R\cos\beta},&\frac{\cos\alpha}{R\cos\beta},& 0\\
-\frac{\cos\alpha}{R}\frac{ \beta+\sin\beta\cos \beta}{\cos\beta+\beta
\sin\beta},&-\frac{\sin\alpha}{R}\frac{ \beta+\sin\beta\cos \beta}{\cos\beta+\beta
\sin\beta},&\frac{  \cos \beta}{2R^2(\cos\beta+\beta
\sin\beta)}
\end{array}\right]
$$

A long, but elementary calculation shows now that, in the new coordinates,
\ieqa\label{deriX}
X=\frac{\cos\beta\cdot\cos(\beta-\alpha)}{\cos\beta+\beta\sin\beta}\partial_R
-\frac{\sin\alpha}{R\cos\beta}\partial_\alpha
-\frac{\beta\cos\alpha+\cos\beta\cdot\sin(\beta-\alpha)}{R(\cos\beta+\beta\sin\beta)}\partial_\beta
\eeqa
and
\ieqa\label{deriY}
Y=-\frac{\cos\beta\cdot\sin(\beta-\alpha)}{\cos\beta+\beta\sin\beta}\partial_R
+\frac{\cos\alpha}{R\cos\beta}\partial_\alpha
-\frac{\beta\sin\alpha+\cos\beta\cdot\cos(\beta-\alpha)}{R(\cos\beta+\beta\sin\beta)}\partial_\beta.
\eeqa
Hence,
\begin{lemma}
\label{lemmaxdist}
Let $\Pi$ be the plane $\{t=0\}$. Then
$$
X\delta_\Pi=-\sin(\beta-\alpha)
$$
and
$$
Y\delta_\Pi=-\cos(\beta-\alpha).
$$
\end{lemma}
\begin{proof} Recalling (\ref{deriX}) and (\ref{deriY}) we get
$X\delta_\Pi=-\sin(\beta-\alpha)$ and $Y\delta_\Pi=-\cos(\beta-\alpha).$
\end{proof}
\begin{remark}
After the computation of $X\delta_\Pi$, we can deduce the value of $Y\delta_\Pi$
keeping in mind (\ref{eqgradrot}).
\end{remark}
\begin{theorem}
\label{Txxdist}
Let $\Pi$ be the plane $\{t=0\}$. Then
$$
\mbox{Hess}_\Hi
\delta_\Pi=\frac{1}{R\cos\beta(\cos\beta+\beta\sin\beta)}\left[\begin{array}
{lll}
 \cos(\alpha-\beta) A,&
 \cos(\alpha-\beta) B\\
 \sin(\alpha-\beta) A,& \sin(\alpha-\beta) B
\end{array}\right],
$$
where
$$
A=-(\sin\alpha\cos\beta+\beta\cos(\alpha+\beta)
+\cos^2\beta\sin(\alpha-\beta))
$$
and

$$
B=(\cos\alpha\cos\beta+
\beta\sin(\alpha+\beta)
+\cos^2\beta\cos(\alpha-\beta)).
$$
Moreover for $P\in \Pi$, $P\not =O$, then
$$
-(\Delta \delta_\Pi)_{\sigma=0}=0,
$$
i.e. the mean curvature is $0$
\end{theorem}

\begin{proof} Recalling Lemma \ref{lemmaxdist} we get by straightforward
computation:
$$
XX\delta_\Pi=-\frac{\sin\alpha}{R\cos\beta}\cos(\alpha-\beta)+
\frac{\beta\cos\alpha-\cos\beta\sin(\alpha -\beta)}
{R(\cos\beta+\beta\sin\beta)}\cos(\alpha-\beta)=
\frac{\cos(\alpha-\beta)}{R\cos\beta}A
$$

$$
YX\delta_\Pi=\frac{\cos\alpha}{R\cos\beta}\cos(\alpha-\beta)+
\frac{\beta\sin\alpha+\cos\beta\cos(\alpha-\beta)}
{R(\cos\beta+\beta\sin\beta)}\cos(\alpha-\beta)=
\frac{\cos(\alpha-\beta)}{R\cos\beta}B
$$
$$
XY\delta_\Pi=-\frac{\sin\alpha}{R\cos\beta}\sin(\alpha-\beta)+
\frac{\beta\cos\alpha-\cos\beta\sin(\alpha-\beta)}
{R(\cos\beta+\beta\sin\beta)}\sin(\alpha-\beta)=
\frac{\sin(\alpha-\beta)}{R\cos\beta}A
$$

$$
YY\delta_\Pi\frac{\cos\alpha}{R\cos\beta}\sin(\alpha-\beta)+
\frac{\beta\sin\alpha+\cos\beta\cos(\beta-\alpha)}
{R(\cos\beta+\beta\sin\beta)}\sin(\alpha-\beta)=\frac{\sin(\alpha-\beta)}{R\cos\beta}B.
$$
In particular, we have
$$\mbox{Trace}(\mbox{Hess}_\Hi
\delta_\Pi)_{\sigma=0}=\mbox{Trace}\left(\frac{1}{R
 }\left[\begin{array} {lll}
 -2\cos\alpha\sin\alpha,&
 2\cos^2\alpha\\
 -2\sin^2\alpha,&2\cos\alpha\sin\alpha
\end{array}\right]\right)=0.
$$
As a consequence the mean curvature of $\Pi-O$ is $0$. This fact could also be
proved as follows. The punctured plane $\Pi-O$ is the union of
straight half line leaving $O$, hence its Heisenberg curvature at a
point $P$ is the Euclidean curvature of a straight line. by Theorem \ref{Tcurveuc}.
\end{proof}
\begin{remark}
Notice that $Hess_\Hi
\delta_\Pi$ is singular as $R\to 0$  and as $\beta\to\pm\pi/2$. For
instance, if $\alpha=0$, then
$YY\delta_\Pi$ diverges as $\beta\to\pm\pi/2$ and $R\to 0$. Indeed
\ieqas
&&(\mbox{Hess}_\Hi
\delta_\Pi)_{\alpha=0}\\
&=&\frac{1}{R\cos\beta(\cos\beta+\beta\sin\beta)}
\left[\begin{array}
{lll}
 -(\beta -\cos\beta\sin\beta)\cos^2 \beta ,&
  (\cos\beta +\beta\sin\beta+\cos^3\beta)\cos\beta\\
 (\beta -\cos\beta\sin\beta)\sin\beta\cos\beta,&
 -(\cos\beta +\beta\sin\beta+\cos^3\beta)\sin\beta
\end{array}\right]
\eeqas
and
$$
(YY
\delta_\Pi)_{\alpha=0}=-\frac{\tan \beta}{R}\frac{\cos\beta
+\beta\sin\beta+\cos^3\beta }{\cos\beta+\beta\sin\beta} =-\frac{\tan
\beta}{R}(1+\frac{\cos^3\beta }{\cos\beta+\beta\sin\beta}).
$$
\end{remark}
The singularity as $R\to0$ is expected, by Theorem \ref{Thess}. The
singularity as $\beta\to\pm\pi/2$ shows that the estimate at the end
of Theorem \ref{Thess} can not be extended outside the surface.

\subsection{The Carnot-Charath\'eodory sphere and the Hessian of the
  distance function}
Here we consider $\Sigma=\partial B(0,1)$, the sphere of unitary
radius with respect to the Carnot-Charath\`eodory metric. The poles of
$\Sigma$ are $NP=(0,0,1/\pi)$ and $SP=(0,0,1/\pi)$. The calculation of
the horizontal Hessian of $\delta_\Sigma$ will be accomplished in several
steps. First, we find the equations of the horizontal curves lying on
$\Sigma$. Via Theorem \ref{Tcurveuc}, this allows us to compute the
mean curvature of $\Sigma$ at every point. By Theorem \ref{Thess},
(\ref{eqhessns}) and (\ref{eqcurvQ}), and Corollary \ref{Tcurv}, the second order terms of the
formula relating $\mbox{Hess}_\Hi$ to the equation of $\Sigma$ are all
contained in the expression for the mean curvature of $\Sigma$. Hence,
in order to compute $\mbox{Hess}_\Hi \delta_\Sigma$ we are left with the
calculation of some first order derivatives.

Finally, we note that the calculation of $\mbox{Hess}_\Hi \delta_\Sigma$
immediately leads to an expression for the horizontal Hessian of the
Carnot-Charath\'eodory distance in $\Hi$, $\mbox{Hess}_\Hi d(O,\cdot)$.

\smallskip

 We can  parametrize $\Sigma$ in several useful ways. Let
 $(z,t)\in\Sigma$. Then, there are $\phi\in[0,2\pi)$ and
 $u\in[-2\pi,2\pi]$ such that
\begin{equation}
\label{eqsigmauphi}
\begin{cases}
z(u)=e^{i\phi}2\frac{\sin(u/2)}{u}\cr
t(u)=2\frac{u-\sin(u)}{u^2}
\end{cases}
\end{equation}
The parametric representation in (\ref{eqsigmauphi}) can be easily
deduced by the geodesic equations (\ref{eqGeo}).

Below we consider the upper half of $\Sigma$, corresponding to $u\ge0$, so that 
$B(O,1)$ is below the surface.

We see, now, how $\Sigma$ is ruled by horizontal
curves. Incidentally, this provides a differently useful
parametrization of $\Sigma$.
\begin{proposition}
\label{proprule}
The surface $\Sigma-\{NP,SP\}$ is the disjoint union of horizontal
curves $\gamma_\theta(u)=(z(u),t(u))$ ($|u|<2\pi$),
\begin{equation}
\label{eqhsigma}
\begin{cases}
z(u)=e^{i(\psi(u/2)+\theta)}2\frac{\sin(u/2)}{u}\cr
t(u)=2\frac{u-\sin(u)}{u^2}
\end{cases}
\end{equation}
where $\psi$ satisfies the conditions $\psi(0)=0$ and
\begin{equation}
\label{eqhprime}
\psi^\prime(p)=\cot^2p-\frac{\cot p}{p}
\end{equation}
These are the only horizontal curves on $\Sigma$.
\end{proposition}
In (\ref{eqhprime}), observe that $\psi^\prime$ is well defined for $p=0$
as well.
\begin{proof} First, note that, by (\ref{eqsigmauphi}), any curve
  parametrized as in (\ref{eqhsigma}) lies on $\Sigma$. Let now
  $\gamma$ be the curve in (\ref{eqhsigma}) corresponding to $\theta=0$
  (by rotation invariance,
  $R_\theta\gamma=\gamma_\theta$ is horizontal if and only if $\gamma$
  is). The curve $\gamma=(x,y,t)$ is horizontal if and only if
\begin{equation}
\label{equgamma}
\dot{t}(u)=2\dot{x}(u)y(u)-2x(u)\dot{y}(u).
\end{equation}
Here,
$$
\begin{cases}
x(u)=\cos(\psi(u/2))\frac{\sin(u/2)}{u/2}\cr
y(u)=\sin(\psi(u/2))\frac{\sin(u/2)}{u/2}
\end{cases}
$$
Let $p=u/2$
We make (\ref{equgamma}) into an explicit O.D.E. A calculation shows
that
$$
2\dot{x}(2p)y(2p)-2x(2p)\dot{y}(2p)=-\psi^\prime(p)\frac{\sin^2(p)}{p^2}
$$
and that
$$
\dot{t}(2p)=\cos(p)\frac{\sin(p)-p\cos(p)}{p^3}
$$
Equation (\ref{eqhprime}) follows immediately.
\end{proof}

It can be proved that $\psi(p)\to\infty$ as $p\to\pi$, i.e. that the
projection of $p(\gamma_\theta)$ onto the plane $\{t=0\}$ spiralizes infinitely
many times around the origin. It can also be shown that the Heisenberg
length of $\gamma_\theta$, i.e. the Euclidean length of $p(\gamma_\theta)$,  is infinite.

Next, we compute the mean curvature of $\Sigma$.
\begin{proposition}
\label{propcurvsigma}
The curvature of $\Sigma$ at $P=(z,t)$, where
$$
\begin{cases}
|z|=2\frac{\sin(u/2)}{u}\cr
t=2\frac{u-\sin(u)}{u^2}
\end{cases}
$$
is the number
\begin{equation}
\label{eqcurvpgamma}
h_\Sigma(P)=\frac{u\cos u-\sin  u}{\frac{u}{2}\cos\left(\frac{u}{2}\right)-\sin\left(\frac{u}{2}\right)}
\frac{u}{4\sin\left(\frac{u}{2}\right)}
\end{equation}
\end{proposition}
\begin{proof} By Theorem \ref{Tcurveuc}, we only need to apply the formula for the computation of the mean
  curvature of a curve in parametric form
(see, e.g., \cite{GT}) to the curves described in Proposition \ref{proprule}.
\end{proof}

By Proposition \ref{propcurvsigma}, in order to compute the Hessian of $\delta_\Sigma$, we are
left with the calculation of a number of {\it first order} derivatives
of $f$.
\begin{lemma}
\label{lemmasigmafirst}
Let $t=f(x,y)$ be the equation of the upper half of $\Sigma$, let
$(z,t)$ be as in (\ref{eqsigmauphi}) and $u=2v$. Then,
\begin{equation}
\label{eqsigmafirst}
\begin{cases}
Xf(z,t)=2\cos\phi\frac{\cos(\theta)}{v}+2\sin\phi\frac{\sin(v)}{v})=\frac{2}{v}\cos(\phi-v)\cr
Yf(z,t)=2\sin\phi\frac{\cos(\theta)}{v}-2\cos\phi\frac{\sin(v)}{v}=\frac{2}{v}\sin(\phi-v)
\end{cases}
\end{equation}
\end{lemma}
\begin{proof}
Let $t=f(x,y)$.
By the chain rule,
\begin{equation}
\label{eqchain}
J\left(\begin{array}{ll}
t&\\
x&y
\end{array}\right)
=J\left(\begin{array}{ll}
t&\\
u&\phi
\end{array}\right)
\cdot
J\left(\begin{array}{ll}
x&y\\
u&\phi
\end{array}\right)^{-1}
\end{equation}
By (\ref{eqsigmauphi}),
$$
J\left(\begin{array}{ll}
t&\\
u&\phi
\end{array}\right)=\left(\frac{\cos(u/2)}{(u/2)^3}\left[\sin(u/2)-u/2\cos(u/2)\right],0\right)
=\left(-\frac{\cos(u/2)}{u/2}h(u/2),0\right)
$$
where
$$
h(s)=\frac{s\cos(s)-\sin(s)}{s^2}
$$
and
$$
J\left(\begin{array}{ll}
x&y\\
u&\phi
\end{array}\right)=\left[\begin{array}{ll}
\frac{1}{2}h(u/2)\cos\phi&-\frac{\sin(u/2)}{u/2}\sin\phi\\
\frac{1}{2}h(u/2)\sin\phi&\frac{\sin(u/2)}{u/2}\cos\phi
\end{array}\right]
$$
so that
$$
J\left(\begin{array}{ll}
x&y\\
u&\phi
\end{array}\right)^{-1}=\frac{1}{\frac{1}{2}h(u/2)\frac{\sin(u/2)}{u/2}}
\left[\begin{array}{ll}
\frac{\sin(u/2)}{u/2}\cos\phi&\frac{\sin(u/2)}{u/2}\sin\phi\\
-\frac{1}{2}h(u/2)\sin\phi&\frac{1}{2}h(u/2)\cos\phi
\end{array}\right]
$$
Using (\ref{eqchain}), we obtain
$$
J\left(\begin{array}{ll}
t&\\
x&y\end{array}\right)=
-\left(\frac{2\cos\phi}{u},\frac{2\sin\phi}{u}\right)
$$
and the conclusion of the lemma follows immediately.
\end{proof}

We have now all the ingredients for the calculation of
$\mbox{Hess}_\Hi \delta_\Sigma$. 
\begin{theorem}
\label{Thesssigma}
Let $P=(z,t)$ be the point of the upper half of $\Sigma$,
\begin{equation}
\label{eqsigmauphii}
\begin{cases}
z(u)=2\frac{\sin(u/2)}{u}\cr
t(u)=2\frac{u-\sin(u)}{u^2}
\end{cases}
\end{equation}
Then,
\begin{equation}
\label{eqhessnsigma}
\mbox{Hess}_\Hi\ \delta_\Sigma(P)=v_P\Sigma\otimes v_P\Sigma \cdot (h_\Sigma(P)I+p_\Sigma(P)J),
\end{equation}
where
$$
v_P\Sigma=\sin(\phi-v)\cdot X_P-\cos(\phi-v)\cdot Y_P,
$$
$$
h_\Sigma(P)=
\frac{u\cos u-\sin  u}{\frac{u}{2}\cos\left(\frac{u}{2}\right)-\sin\left(\frac{u}{2}\right)}
\frac{u}{4\sin\left(\frac{u}{2}\right)},
$$
and
$$
p_\Sigma(P)=2u.
$$
\end{theorem}
\begin{proof}
The theorem follows immediately from Theorem \ref{Tmain}, Proposition \ref{propcurvsigma}, 
Lemma \ref{lemmasigmafirst}, (\ref{reggio}) and the definition of $v_P\Sigma$.
\end{proof}

\bigskip

Nicola Arcozzi

Dipartimento di Matematica dell'Universit\`a

Piazza di Porta S.~Donato, 5, 40126 Bologna, Italy

E-mail: arcozzi@dm.unibo.it

\vspace{1cm}

Fausto Ferrari

Dipartimento di Matematica dell'Universit\`a

Piazza di Porta S.~Donato, 5, 40126 Bologna, Italy

and

C.I.R.A.M.

Via Saragozza, 8, 40123 Bologna, Italy.

E-mail: ferrari@dm.unibo.it

\end{document}